\documentclass{article}

\usepackage[utf8]{inputenc}
\usepackage{amsmath,amssymb,amsfonts}
\usepackage{graphicx}
\usepackage{authblk}
\usepackage{palatino}
\usepackage{appendix}
\usepackage{subcaption}
\usepackage{natbib}

\title{Adaptive Network Modeling of Social Distancing Interventions}
\author[1]{Carl Corcoran \thanks{Corresponding Author. Email: ctcorcoran@ucdavis.edu}}
\author[2]{John Michael Clark}
\affil[1]{Department of Mathematics, University of California, Davis, Davis, CA, USA}
\affil[2]{Department of Mathematics, Oklahoma State University, Stillwater, OK, USA}

\begin{document}
\maketitle

\begin{abstract}
    \noindent The COVID-19 pandemic has proved to be one of the most disruptive public health emergencies in recent memory. Among non-pharmaceutical interventions, social distancing and lockdown measures are some of the most common tools employed by governments around the world to combat the disease. While mathematical models of COVID-19 are ubiquitous, few have leveraged network theory in a general way to explain the mechanics of social distancing. In this paper, we build on existing network models for heterogeneous, clustered networks with random link activation/deletion dynamics to put forth realistic mechanisms of social distancing using piecewise constant activation/deletion rates. We find our models are capable of rich qualitative behavior, and offer meaningful insight with relatively few intervention parameters. In particular, we find that the severity of social distancing interventions and when they begin have more impact than how long it takes for the interventions to take full effect.
\end{abstract}

\section{Introduction} \label{intro}

The global COVID-19 pandemic has upended modern life an placed an enormous epidemiological, economic, and social burden on the world's resources. The gravity of events has brought the need for epidemiological modeling into sharp focus. As the pandemic spread around the world in the absence of a vaccine, non-pharmaceutical interventions including social distancing, quarantine, and lockdown measures proliferated, and bringing these interventions into modeling efforts has remained paramount. 

\begin{figure}[h]
    \centering
    \begin{subfigure}{.3\textwidth}
        \centering
        \includegraphics[width=.8\textwidth]{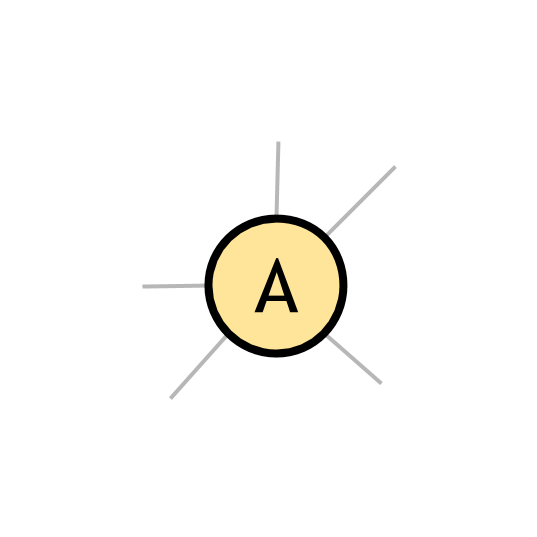}
        \caption{}
    \end{subfigure}
    \begin{subfigure}{.3\textwidth}
        \centering
        \includegraphics[width=.8\textwidth]{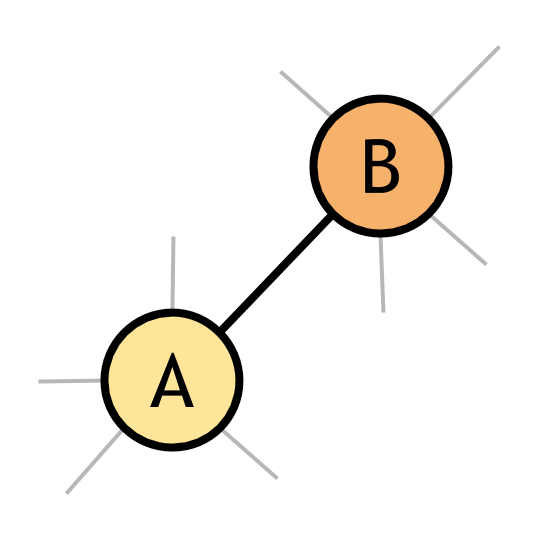}
        \caption{}
    \end{subfigure}
    \begin{subfigure}{.3\textwidth}
        \centering
        \includegraphics[width=\textwidth]{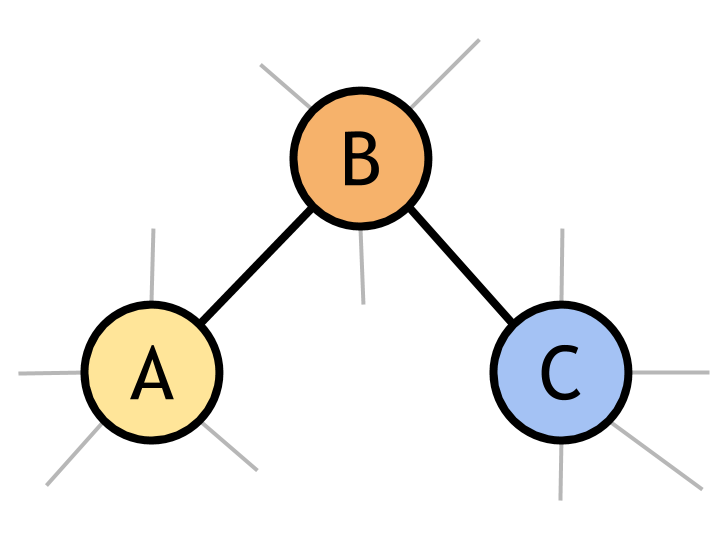}
        \caption{}
    \end{subfigure}
    \caption{Diagrams of network structures whose evolution is modeled by the pairwise model: (a) node in state $A$, (b) pair in state $A-B$, (c) triple in state $A-B-C$.}
    \label{fig:nodepairtriple}
\end{figure}

In recent years, network-based models of epidemic spread have become an increasingly popular paradigm \citep{pastor-satorras_epidemic_2015,kiss_mathematics_2017}, and network science generally has been recognized for its potential to contribute solutions to the current crisis \citep{eubank_commentary_2020}. Most network models represent individuals as nodes in a network, and their contacts as edges connecting the nodes. Moreover, many models assume that the network is static--that the edges between nodes don't change over time--and thus the epidemic spreads from node to node across these edges. Among static network models, pairwise models \citep{keeling_effects_1999-1,eames_modeling_2002} are both frequently used and well-studied. Pairwise models track not only the number of nodes in a given state, but pairs, triples, and higher order motifs as well (Fig. \ref{fig:nodepairtriple}). An advantage of pairwise models in that in their full form, they exactly model (in expectation) the continuous time Markov chain formulation of epidemic spread on a network \citep{taylor_markovian_2012}. 

Pairwise models have been successfully applied to a number of disease natural histories and different network types. Two important network features that play a role in the theory of pairwise models are degree heterogeneity and clustering. The degree of a node in a network is the number of edges to which it is connected, and the degree distribution is the probability distribution of selecting a random node with a given degree. The degree distribution plays a fundamental role in many network models, and is particularly powerful when described as a probability generating function. The clustering coefficient is the ratio of triangles to connected triples in the network. While clustering in an important component of network structure, it has not widely been incorporated to pairwise models. We acknowledge two major benefits of degree heterogeneous, clustered models.
First, including both or either as modeling consideration affects epidemic dynamics in a nontrivial way \citep{house_insights_2011,keeling_effects_1999-1} and second, both have been shown to be features of realistic contact networks \citep{read_dynamic_2008}.

Though static networks model some forms of complexity well, an important aspect of real-world contact networks is that some connections change in response to disease dynamics or public health measures. By relaxing the static network assumption, dynamic or ``adaptive" network models \citep{gross_adaptive_2009} can capture both the dynamics of the network and the epidemic dynamics on the network. A number of models have been recently proposed for that describe a variety of network dynamic processes. \citet{gross_epidemic_2006} introduced a model of edge rewiring, where susceptible individuals break connections with infectious individuals and reconnect to susceptible individuals at random. Another model for network dynamics is random link addition/deletion \citep{kiss_modelling_2012} where individuals break and form new contacts and constant rates. Their approach is notable for its intuitiveness as a simple dynamic model, and also its use of probability generating functions as a tool to describe network dynamics. A related model is link addition/deletion on a fixed network \citep{tunc_epidemics_2013,shkarayev_epidemics_2014}, where individuals can temporarily deactivate contacts with infectious individuals, and reactivate them when their contact is not infectious. While much of the focus of the adaptive network literature has been involved in analyzing the resulting dynamical systems, particularly for SIS-type diseases, some works have focused on the role of network dynamics in controlling or mitigating epidemic spread \citep{youssef_mitigation_2013,selley_dynamic_2015}.

Network models in general offer a compromise between two other common modeling techniques: compartment models and agent-based simulations. They are able to capture more complex contact structure than simple compartment models, while offering analytical tractability that many agent-based simulations lack. Despite this, models of non-pharmaceutical interventions have tended to favor simulation or compartment models \citep{ahmed_effectiveness_2018,davey_effective_2008}. In the early stages of the COVID-19 pandemic, complex individual-based simulations offered major insights about the effectiveness of non-pharmaceutical interventions \citep{ferguson_report_2020}. However, the high computational cost can make investigating the impacts of intervention policies with a large number of parameters a challenging endeavor. Network models, especially those with a relatively small number of equations, can offer broad insights at reduced cost. While some models of social distancing have incorporated contact network structure as a major consideration \citep{valdez_intermittent_2012,glass_targeted_2006}, differential equation network models of such interventions are uncommon. Adaptive network models in particular can offer a new perspective on questions surrounding social distancing and other non-pharmaceutical interventions made pressing by the COVID-19 pandemic.

In this paper, we develop simple, novel mechanisms to incorporate social distancing into a network model of epidemic spread, using COVID-19 as the central case study to investigate the impact of a range of interventions. First, we develop a pairwise SEIR model with random link activation/deletion dynamics--that is edges are added and deleted at constant rates independent of the epidemic dynamics on the network. Furthermore, the model incorporates degree heterogeneity and clustering, which offers increased realism over simpler network or compartment models. To apply the model, we use bipartite mixing networks to generate large heterogeneous, clustered contact networks coupled with disease dynamics given by epidemiological parameters estimated for COVID-19. Next, we develop two mechanisms of social distancing using piecewise constant link activation and deletion rates. The first is a single intervention event, where the average number of contacts decreases, is held constant, and then recovers; the second allows for multiple interventions which restart depending on the prevalence of the disease. While we investigate the implications of these policies for COVID-19 on a specific type of heterogeneous, clustered network, both the adaptive network model and the social distancing schemes are more generally applicable to a variety of networks and epidemiological parameters. Finally, we consider the public health implications of the latter model, finding that certain intervention parameters are more important than others in achieving an effective reduction in overall infections.

\section{Model} \label{model}

To begin construction of the full model, we consider SEIR dynamics on a static network. Pairwise equations for an SEIR epidemic can be found in \citet{keeling_correlation_1997} and \citet{mcglade_correlation_1999}. Model variables include the expected number of susceptible, exposed, infectious, and recovered nodes ($[S],[E],[I]$ and $[R]$ respectively) as well as the expected number of pairs in each state. For example, $[SS]$ is the expected number of connected pairs of susceptible nodes, while $[SI]$ is the expected number of connected pairs of susceptible and infectious nodes. The expected number of connected triples is also considered ($[SSI],[ESI],[ISI]$), though differential equations for these variables are not written. The full SEIR pairwise model is
\begin{align}
\dot{[S]} &= -\beta[SI], \label{eq:SEIR1}\\
\dot{[E]} &= \beta[SI]-\eta[E],\\
\dot{[I]} &= \eta[E]-\gamma[I],\\
\dot{[SS]} &= -2\beta[SSI],\\
\dot{[SE]} &= \beta[SSI]-\beta[ESI] - \eta[SE],\\
\dot{[SI]} &= \eta[SE] -\beta[SI] - \beta[ISI]-\gamma[SI],\\
\dot{[EE]} &= 2\beta[ESI]-2\eta[EE],\\
\dot{[EI]} &= \beta[ISI]+\beta[SI]+\eta[EE]-(\gamma+\eta)[EI],\\
\dot{[II]} &= 2\eta[EI]-2\gamma[II],\label{eq:SEIR9}
\end{align}
where $\beta$ is the transmission rate, $\gamma$ is the recovery rate, and $\eta$ is the rate at which exposed individuals become infectious. The nodes and edges also obey conservation equations
\begin{equation}
N = [S]+[E]+[I]+[R] 
\end{equation}
and
\begin{align}
\langle k \rangle N &= [SS] + [EE] + [II] + [RR] \nonumber \\
&\qquad + 2\left([SE]+[SI]+[SR]+[EI]+[ER]+[IR]\right) 
\end{align}
where $N$ is the number of nodes and $\langle k \rangle$ is the average degree of the network. We note that with the conservation equations, we do not need terms of the type $[AR]$ to determine the evolution of $[S],[E],[I],$ and $[R].$ 

The full model requires dynamical equations for triples of the form $[ASI]$ and higher order motifs as well, leading to a system that is prohibitively large for computations. To make the model tractable, we approximate the expected number of triples $[ASI]$ in terms of pairs and individual nodes, thus closing the system (\ref{eq:SEIR1})-(\ref{eq:SEIR9}). An approximation of this kind is referred to as a triple closure. For triples of the type $A-S-I,$ \citet{house_insights_2011} give a triple closure approximation as
\begin{equation}\label{eq:triple_approx}
[ASI] \approx [AS][SI]\frac{\sum_k (k^2-k)[S_k]}{\left(\sum_k k[S_k]\right)^2}\left(1-\phi \\
+ \phi \frac{N \langle k \rangle [AI]}{\left(\sum_k k[A_k]\right)\left(\sum_k k[I_k]\right)}\right)
\end{equation}
where $[A_k]$ is the expected number of nodes in state $A$ with degree $k.$ Using the network degree distribution probability generating function and introducing new dynamical variables, they develop an SIR model for heterogeneous, clustered networks. In Appendix \ref{AppA}, we derive an analogous heterogeneous, clustered SEIR model complete with link activation and deletion. However, the model complexity induced by (\ref{eq:triple_approx}) is not necessary to accurately capture the combined epidemic and network dynamics, and thus a simpler triple closure will suffice.

A simple yet useful assumption is that degree and state are independent, and thus $[A_k] = p_k [A]$ where $p_k$ is the proportion of nodes with degree $k.$ With this assumption, the resulting triple closure becomes:
\begin{equation}
[ASI] \approx \frac{\langle k^2 - k \rangle}{\langle k \rangle^2}\frac{[AS][SI]}{[S]}\left(1-\phi + \\
\phi \frac{N}{\langle k \rangle}\frac{[AI]}{[A][I]}\right),
\end{equation}
where $\langle k \rangle = \sum_{k=0}^{N-1}kp_k$ and $\langle k^2-k \rangle = \sum_{k=0}^{N-1}(k^2-k)p_k$. We note that if a homogeneous degree distribution is assumed, the closure reduces to clustered closure from \citet{keeling_effects_1999-1}.

With the static model closed, we now incorporate the effects of network dynamics. \citet{kiss_modelling_2012} introduced a simple model of network dynamics, termed random link activation/deletion (RLAD). In this model, independent of epidemic dynamics nonexistent edges are added to the network (or activated) at a constant rate $\alpha$ and existing edges are removed from the network (or deleted) at a constant rate $\omega$. Ignoring epidemic spread and looking at the effects of activation/deletion only, the equation for edges of type $[AA]$ is

\begin{equation} \label{eq:[AA]}
    [\dot{AA}] = \alpha\left([A]([A]-1)-[AA]\right) - \omega[AA]
\end{equation}
and for type $[AB]$ we have
\begin{equation} \label{eq:[AB]}
    [\dot{AB}]=\alpha\left([A][B]-[AB]\right)-\omega[AB].
\end{equation}

Next, we have to consider the effect of activation/deletion on the now time-dependent network quantities: degree distribution moment terms $\langle k \rangle(t), \langle k^2-k\rangle(t)$ and the clustering coefficient $\phi(t).$  Following the example of \citet{kiss_modelling_2012}, dynamical equations for the first two can be easily derived by finding the partial differential equation for the degree distribution generating function 
\begin{equation}
    g(x,t) = \sum_{k=0}^{N-1}p_k(t)x^k.
\end{equation}
The Kolmogorov equations describe the evolution of $p_k(t),$ the proportion of degree $k$ nodes at time $t:$

\begin{equation} \label{eq:kolmogorov}
    \dot{p_k} = \alpha\left(N-k\right)p_{k-1} - \left(\alpha(N-1-k)+\omega k\right)p_k +\omega(k+1)p_{k+1}.
\end{equation}
With some straightforward algebra, we derive a partial differential equation for the degree distribution generating function:
\begin{equation} \label{eq:pde}
    \frac{\partial g}{\partial t} = (x-1)\left(\alpha(N-1)g-(\alpha x+\omega)\frac{\partial g}{\partial x} \right).
\end{equation}
The network quantities $\langle k \rangle$ and $\langle k^2-k\rangle$ can be computed from the generating function as $\langle k \rangle = g_x(1,t)$ and $\langle k^2-k\rangle=g_{xx}(1,t)$. Then, from (\ref{eq:pde}) we derive the dynamical equations:
\begin{align}
    \langle \dot{k}\rangle &= \alpha(N-1)-(\alpha + \omega)\langle k \rangle,\\
    \langle \dot{k^2-k}\rangle &= 2\alpha(N-2)\langle k \rangle - 2(\alpha + \omega)\langle k^2-k\rangle.
\end{align}
The clustering coefficient is defined as the ratio of triangles to connected triples in the network. To compute $\dot{\phi}$, we start with the Kolmogorov equations for $q_k(t),$ the probability that there are $k$ triangles in the network at time $t:$
\begin{equation}
    \dot{q_k} = \alpha(L-3(k-1))q_{k-1}-(\alpha(L-3k)+3\omega k)q_k+3\omega(k+1)q_{k+1}
\end{equation}
where $L = N\langle k^2-k\rangle/2$ is the number of connected triples. From this we derive the differential equation for the expected number of triangles $\langle T \rangle$ as
\begin{equation}
\langle \dot{T}\rangle = \alpha L-3(\alpha+\omega)\langle T\rangle,
\end{equation}
and compute the equation for the clustering coefficient $\phi(t):$ 
\begin{equation}\label{eq:phi}
    \dot{\phi} = 3\alpha - \left(\alpha + \omega+2\alpha(N-2)\frac{\langle k \rangle}{\langle k^2-k\rangle}\right)\phi.
\end{equation}
Finally, we have a full set of equations for a pairwise SEIR for a heterogeneous, clustered network with random link activation and deletion:
\begin{align}
\dot{[S]} &= -\beta[SI], \label{eq:mod1}\\
\dot{[E]} &= \beta[SI]-\eta[E],\\
\dot{[I]} &= \eta[E]-\gamma[I],\\
\dot{[SS]} &= -2\beta[SSI]+\alpha[S]([S]-1)-(\alpha+\omega)[SS], \\
\dot{[SE]} &= \beta[SSI]-\beta[ESI] - \eta[SE]+\alpha[S][E]-(\alpha+\omega)[SE],\\
\dot{[SI]} &= \eta[SE] -\beta[SI] - \beta[ISI]-\gamma[SI]+\alpha[S][I]-(\alpha+\omega)[SI],\\
\dot{[EE]} &= 2\beta[ESI]-2\eta[EE]+\alpha[E]([E]-1)-(\alpha+\omega)[EE],\\
\dot{[EI]} &= \beta[ISI]+\beta[SI]+\eta[EE]-(\gamma+\eta)[EI]\nonumber\\
&\qquad+\alpha[E][S]-(\alpha+\omega)[EI],\\
\dot{[II]} &= 2\eta[EI]-2\gamma[II]+\alpha[I]([I]-1)-(\alpha+\omega)[II],\\
\langle \dot{k}\rangle &= \alpha(N-1)-(\alpha + \omega)\langle k \rangle, \label{eq:k}\\
\langle \dot{k^2-k}\rangle &= 2\alpha(N-2)\langle k \rangle - 2(\alpha + \omega)\langle k^2-k\rangle,\\
\dot{\phi} &= 3\alpha - \left(\alpha + \omega+2\alpha(N-2)\frac{\langle k\rangle}{\langle k^2-k\rangle}\right)\phi, \label{eq:mod12}
\end{align}
where 
\begin{align}
    [SSI] &= \frac{\langle k^2 - k \rangle}{\langle k \rangle^2}\frac{[SS][SI]}{[S]}\left(1-\phi + \phi \frac{N}{\langle k \rangle}\frac{[SI]}{[S][I]}\right),\\
    [ESI] &= \frac{\langle k^2 - k \rangle}{\langle k \rangle^2}\frac{[SE][SI]}{[S]}\left(1-\phi + \phi \frac{N}{\langle k \rangle}\frac{[EI]}{[E][I]}\right),\\
    [ISI] &= \frac{\langle k^2 - k \rangle}{\langle k \rangle^2}\frac{[SI]^2}{[S]}\left(1-\phi + \phi \frac{N}{\langle k \rangle}\frac{[II]}{[I]^2}\right).
\end{align}
\begin{figure}[h]
    \begin{subfigure}[b]{0.49\textwidth}
    \centering
    \includegraphics[width=\textwidth]{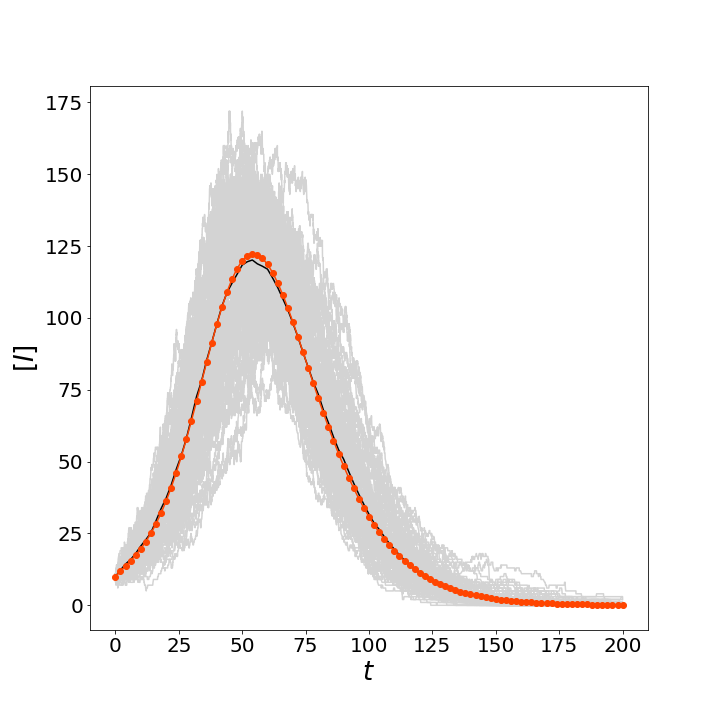}
    \caption{}
    \end{subfigure}
    \begin{subfigure}[b]{0.49\textwidth}
    \centering
    \includegraphics[width=\textwidth]{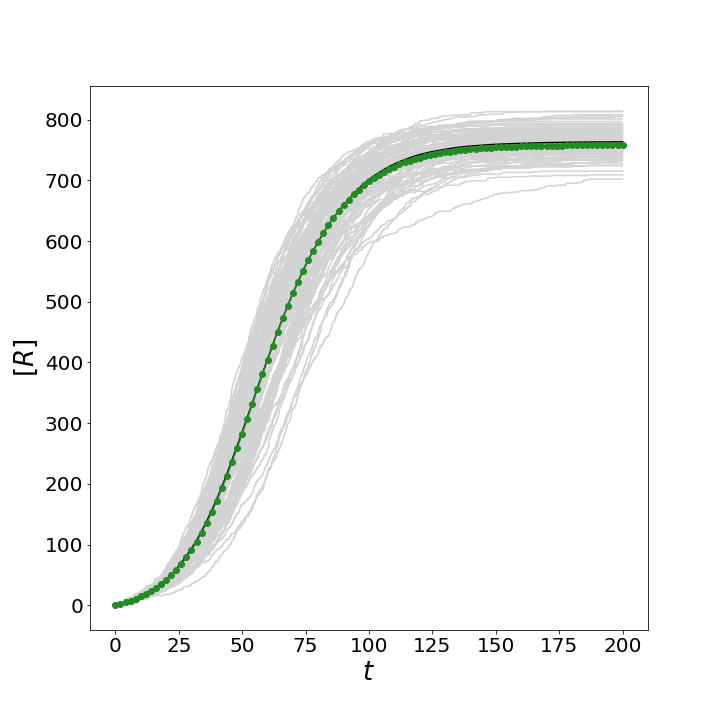}
    \caption{}
    \end{subfigure}
    \caption{Comparison of the model to simulation. 100 trials were run on a unipartite contact network generated from a bipartite network with Poisson degree distributions and $N=500,M=125,\lambda = 4$. Initial conditions are $[E]_0=[I]_0=10,[S]_0=N-20,[R]_0=0.$ Epidemiological and network  parameters are $R_0 = 2.4,\eta = 1/5,\gamma = 1/10,\alpha \approx 2.3\times10^{-5},\omega = 3.4\times 10^{-5}$. Individual simulations are shown in light gray with the mean in black. Model results are (a) $[I](t)$, red circles, (b) $[R](t)$, green circles.}
    \label{fig:sim}
\end{figure}
To demonstrate the validity of this model, we test it against numerical simulations (Fig. \ref{fig:sim}) on a heterogeneous, clustered network---the construction of which is described in Section \ref{network}. Clearly, the model (\ref{eq:mod1})-(\ref{eq:mod12}) is in excellent agreement with the simulations.

\subsection{Network and Epidemiological Parameters}\label{network}

\begin{figure}
    \centering
    \begin{subfigure}[b]{\textwidth}
    \centering
    \includegraphics[width=\textwidth]{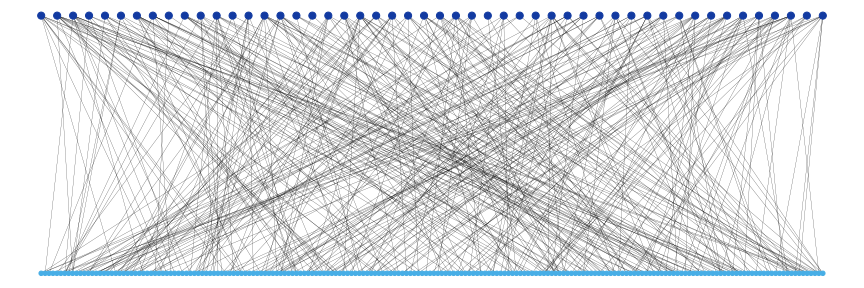}
    \caption{}
    \label{fig:bipartitenetwork}
    \end{subfigure}
    \begin{subfigure}[b]{0.49\textwidth}
    \centering
    \includegraphics[width=\textwidth]{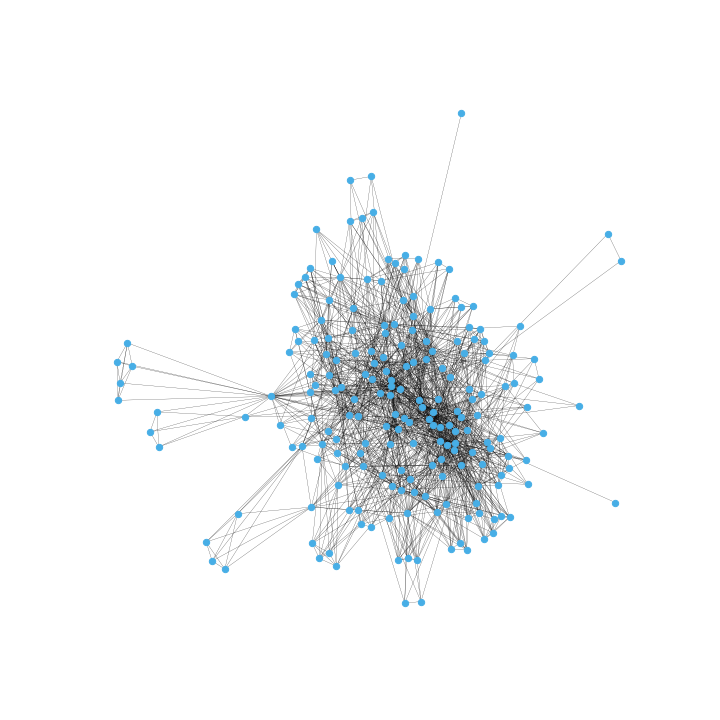}
    \caption{}
    \label{fig:contactnetwork}
    \end{subfigure}
    \centering
    \begin{subfigure}[b]{0.49\textwidth}
    \includegraphics[width=\textwidth]{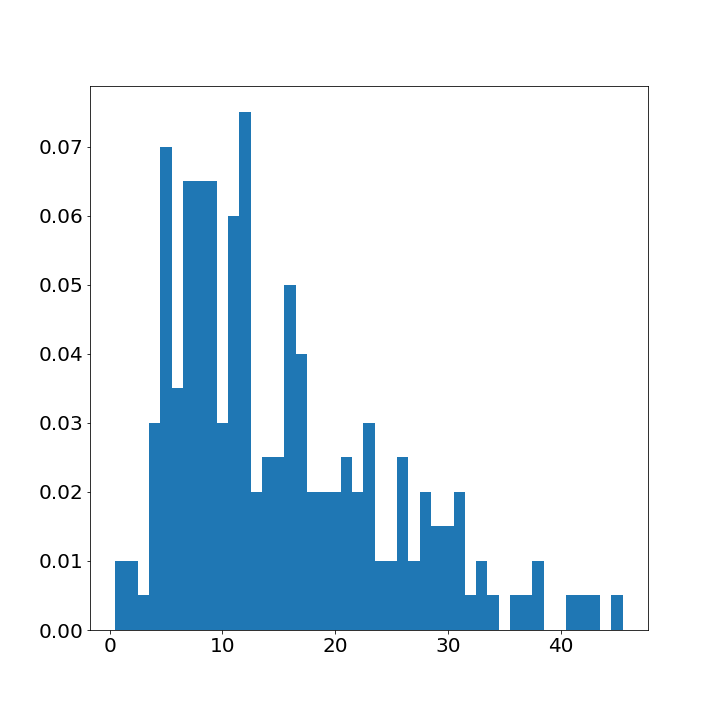}
    \caption{}
    \label{fig:degreedist}
    \end{subfigure}
    \caption{Example contact network (b) and its degree distribution (c) generated from a bipartite mixing network (a). Degree distributions for the individuals and mixing locations are Poisson (as described in Section \ref{network}) with $N=200,M=50$, and $\lambda=2$.}
    \label{fig:network}
\end{figure}
The goal of this paper is to investigate social distancing policies through random link activation/deletion dynamics, which are controlled by the activation and deletion rates $\alpha$ and $\omega.$ Moreover, in building intervention schemes in Section \ref{analysis} new parameters are introduced. In order to consistently compare the efficacy of intervention schemes, network and epidemiological parameters are held the same across schemes. As such, we restrict our attention to a particular heterogeneous, clustered network and epidemiological parameters that are plausible for COVID-19. For completeness, other network types and epidemiological parameters are considered in the Appendix \ref{AppB}.

A consistent challenge of network models is constructing realistic contact networks. In particular, degree heterogeneity and significant clustering are observed in real world social networks \citep{read_dynamic_2008}. To construct such a contact network, we consider a bipartite mixing network \citep{eubank_modelling_2004} with $N$ individuals and $M$ mixing locations (Fig. \ref{fig:bipartitenetwork}). Two individuals are in contact if they both connect to the same mixing location, so we form a contact network as the unipartite projection of the bipartite mixing network (Fig. \ref{fig:contactnetwork}). To introduce degree heterogeneity, we construct a bipartite mixing network where both individuals and mixing locations have Poisson degree distributions. The average individual degree $\lambda$ and average mixing location degree $\mu$ are related by 
\begin{equation}
    N\lambda = M\mu,
\end{equation}
so only $N,M,$ and $\lambda$ are needed to characterize this network. Using generating function techniques \citep{newman_random_2001}, we compute

\begin{align}
    \langle k\rangle &= \frac{N}{M}\lambda^2,\\
    \langle k^2-k\rangle &= \left(\frac{N}{
    M}\right)^2\lambda^3(\lambda+1),\\
    \phi &= \frac{1}{\lambda+1},
\end{align}
for the unipartite contact network, which exhibits both degree heterogeneity (Fig. \ref{fig:degreedist}) and clustering. Unless otherwise specified, the networks in this article are generated using $N=10,000, M= 2,500,$ and $\lambda = 4.$ We acknowledge that though we use a bipartite mixing network to generate a heterogeneous, clustered unipartite contact network, our network dynamics are limited to the contact network. Mobility networks \citep{chang_mobility_2021} have been used to great effect for COVID-19, and suggest a fruitful path forward for bipartite network dynamics.

Numerous recent studies have estimated important epidemiological quantities for the spread of Sars-CoV-2, including the length of the incubation period, the length of the infectious period, and the basic reproduction number $R_0$. We choose the plausible estimates in line with recent studies: average incubation period of 5 days \citep{linton_incubation_2020,zhang_evolving_2020}, average infectious period of 10 days \citep{you_estimation_2020}, and $R_0 = 2.4$ \citep{li_early_2020,anastassopoulou_data-based_2020}. To incorporate these into the model, we note that $1/\eta$ and $1/\gamma$ are the average lengths of the incubation and infectious periods respectively, and thus $\eta = 0.2,\gamma = 0.1.$ We do not derive $R_0$ for the model (\ref{eq:mod1})-(\ref{eq:mod12}), but instead consider the basic reproduction number for a heterogeneous, clustered population from \citet{miller_spread_2009}, which is given as the series
\begin{equation}
    R_0 = \frac{\langle k^2-k\rangle}{\langle k \rangle}\frac{\beta}{\beta+\gamma}-\phi\frac{\langle k^2-k\rangle}{\langle k \rangle}\left(\frac{\beta}{\beta+\gamma}\right)^2+\dots \label{eq:R_0}
\end{equation}
Ignoring higher order terms, we can compute $\beta$ from (\ref{eq:R_0}) when $R_0 = 2.5$ With these parameters, we plausibly model the spread of COVID-19 through a moderately sized heterogeneous, clustered population in the following sections, while introducing various social distancing interventions to mitigate or control the epidemic.

\section{Analysis} \label{analysis}

Social distancing and lockdown measures have been used to curb the spread of infectious diseases throughout history, and are some of the most ubiquitous non-pharmaceutical interventions in the current COVID-19 pandemic. Many compartment-based models that incorporate social distancing do so in a phenomenological manner through the transmission rate, but adaptive network models present an opportunity to describe a social distancing mechanism in a fundamental way.
\begin{figure}
    \centering
    \includegraphics[width=\textwidth]{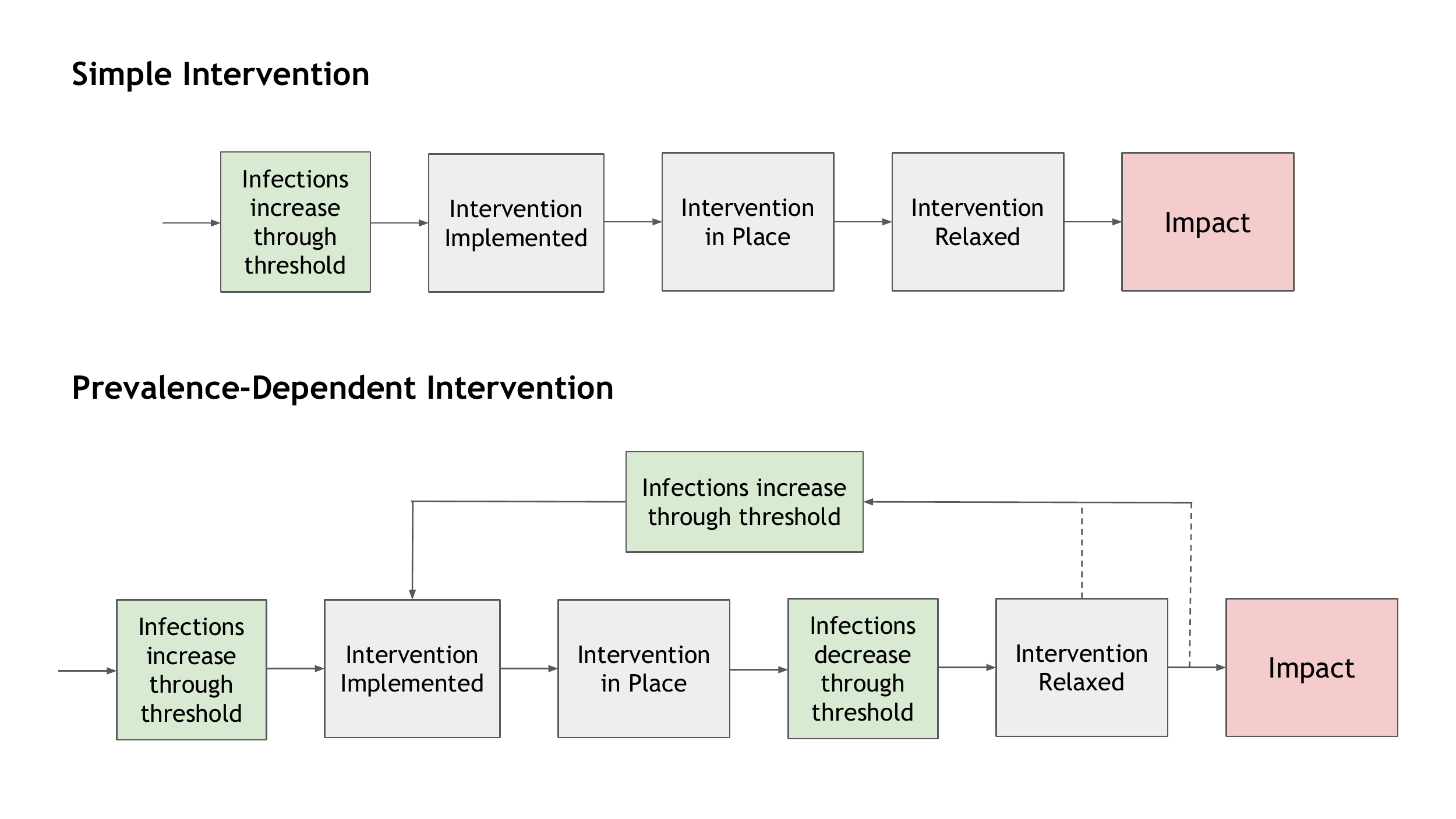}
    \caption{Schematic of the Simple and Prevalence-Dependent Interventions. Both interventions are triggered by a threshold condition, and proceed through the described intervention until the epidemic ends and the impacts of the interventions can be evaluated.}
    \label{fig:flowchart}
\end{figure}
A simple model of such interventions can be naturally characterized by the link activation/deletion process. During periods of social distancing and lockdown, individual contacts break; during periods of relaxation of the measures, individual contacts form. In this section, we develop two social distancing schemes (Fig. \ref{fig:flowchart}). Both social distancing schemes begin when the prevalence $[I](t)$ reaches some specified threshold level. For the simple intervention scheme, contacts break as the intervention is implemented, then contacts stay fixed as the intervention is in place, and finally contacts form until they reach their pre-intervention levels. The prevalence-dependent scheme unfolds similarly, but with two notable differences. First, after the intervention, contacts do not start forming again until the prevalence has dropped below the threshold. Second, any time the prevalence reaches the threshold again, the intervention restarts. This allows for multiple implementations of a social distancing intervention throughout the course of the epidemic. 

Critically, we do not treat these schemes as a mere modeling exercise, but are interested in the impact of each intervention scheme at the end of the epidemic. We develop two simple metrics to evaluate the effectiveness of the simple and prevalence-dependent interventions. First, we consider each intervention's ability to reduce the cumulative number of infections, known as the final size of the epidemic. Second, we also consider how many infections occur above the threshold value for prevalence. These two measures reflect two different yet crucial public health goals, and do not necessarily agree on which interventions are the most effective. Both must be considered to get a complete picture of an intervention's impact. In this section, we derive these two metrics mathematically, and describe the simple and prevalence-dependent interventions while assessing their overall effects.

\subsection{Evaluation Metrics} \label{eval}

The first measure of intervention effectiveness we introduce is the \textit{Relative Change in the Final Size} (RCFS). The final size of an epidemic is the cumulative number of infections that occur over the course of the epidemic. In terms of the model, the final size can be found as the limiting value of the recovered individuals $[R]$: 
\begin{equation*}
    \lim_{t\to\infty}[R](t) = R_{\infty}.
\end{equation*}
We compare the final size of the epidemic with no intervention $R_{\infty}$ to the final size where an intervention has been implemented $R_{\infty}^{\text{int}}.$ We then define the RCFS as
\begin{equation}
\text{RCFS} = \frac{R_{\infty}^{\text{int}} - R_{\infty} }{R_{\infty}}.
\end{equation}
An effective intervention will lead to a decrease in final size, so an RCFS near 0 is unsuccessful, while an RCFS near $-1$ is extraordinarily successful. However, it is important to note that for brief, intense intervention schemes, it is possible that the final size actually increases. In this case, the network parameters change quickly, before significant disease spread, so the epidemic unfolds on a fundamentally different static network.

While the relative change in the final size provides an overall measure of the effectiveness of interventions, reducing cumulative infections alone is not the only public health goal that an intervention scheme might seek to accomplish. In some schemes, a large number of infections occur above the threshold despite a large reduction in the final size of the epidemic. This can be particularly pernicious if the threshold represents some fixed resource such as healthcare capacity, where a large number of infections above the threshold could lead to higher mortality and other negative outcomes. To account for this, we compute the \emph{Cumulative Infections Above Threshold} (CIAT). Let $t_1,t_2,...$ be the sequence of times when $[I] = qN.$ Assuming $[\dot{I}]\neq 0$ at any time in the sequence, the continuity of $[I](t)$ implies that the prevalence is above the threshold on the intervals $[t_{2i-1},t_{2i}]$ for $i = 1,2,3,...$ Thus, the CIAT may be defined as 
\begin{equation}
    \text{CIAT} = \sum_{i}\int_{t_{2i-1}}^{t_{2i}}[I](t) - qN dt.
\end{equation}
We note that the units of CIAT are person-time--for a metric with units of population, we compute the \emph{Average Infections Above Threshold} (AIAT):
\begin{equation}
    \text{AIAT} = \frac{\text{CIAT}}{\sum_{i}t_{2i}-t_{2i-1}} \label{eq:AIAT}
\end{equation}
Using the relation $[\dot{R}] = \gamma[I],$ equation (\ref{eq:AIAT}) becomes
\begin{equation}
    \text{AIAT} = \frac{\sum_{i}[R](t_{2i})-[R](t_{2i-1})}{\gamma\sum_{i}t_{2i}-t_{2i-1}} -qN,
\end{equation}
which is convenient for computations.

\subsection{Simple Intervention} \label{simple}

For a simple model of social distancing, we consider a scheme that unfolds in three successive phases, each with variable length. The effects of the intervention scheme on the contact network are characterized through the average number of contacts $\langle k \rangle (t).$ The intensity of the intervention can be thought of as how severely the average number of contacts are reduced, so we introduce a severity parameter $p \in [0,1).$ 
\begin{figure}[h]
    \centering
    \includegraphics[width=\textwidth]{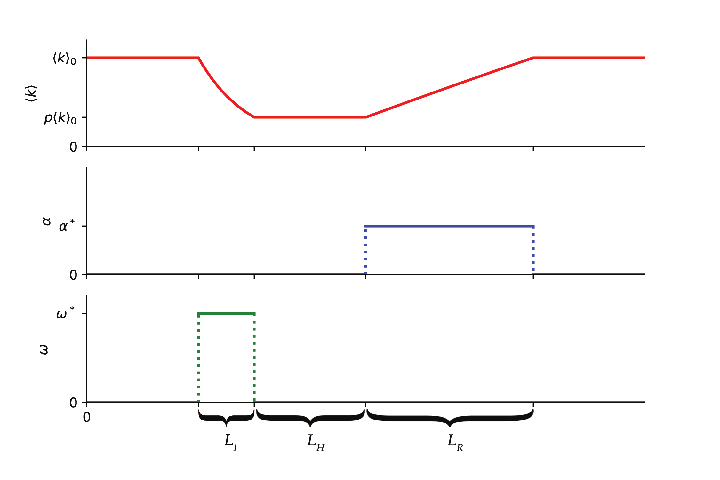}
    \caption{Simple Intervention. Once the intervention begins, edges are deleted at rate $\omega^*$ for $L_I$ days until the average number of contacts $\langle k\rangle$ drops to $p\langle k \rangle_0.$ For the next $L_H$ days, no changes are made to the nework. Then, edges are added at rate $\alpha^*$ for $L_R$ days, until the average number of contacts $\langle k \rangle$ increases to back to $\langle k\rangle_0$}
    \label{fig:intervention1}
\end{figure}
The top panel of Fig. \ref{fig:intervention1} shows how the $\langle k \rangle$ changes over time as the result of the intervention. In the first phase, as social distancing measures are put into place, the average number of contacts decreases from its pre-intervention level $\langle k \rangle_0$ to $p\langle k \rangle_0.$ In the second phase, with the measures fully in place, the average number of contacts remains constant at $p\langle k \rangle.$ In the third phase, social distancing measures are relaxed and the average number of contacts increases to its pre-intervention level $\langle k \rangle_0$. 

\begin{figure}[h]
    \centering
    \begin{subfigure}[b]{.3\textwidth}
    \centering
    \includegraphics[scale=.35]{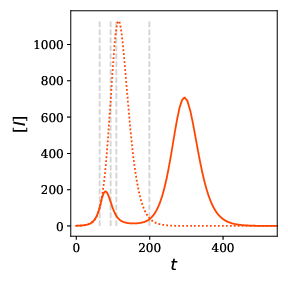}
    \caption{}
    \end{subfigure}
    \begin{subfigure}[b]{.3\textwidth}
    \centering
    \includegraphics[scale=.35]{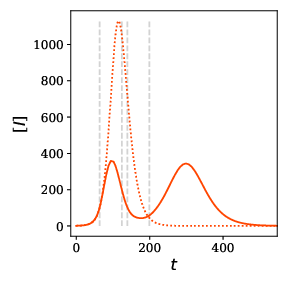}
    \caption{}
    \end{subfigure}
    \begin{subfigure}[b]{0.3\textwidth}
    \centering
    \includegraphics[scale=.35]{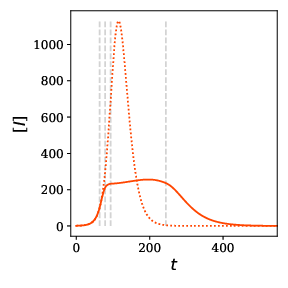}
    \caption{}
    \end{subfigure}
    \caption{Example infection curves $[I](t)$ for the simple intervention with $q=0.01.$ The other intervention parameters are (a) $p = 0.125$, $L_I = 30$, $L_R = 90$, (b) $p = 0.25$, $L_I = 60$, $L_R = 60$, and (c) $p = 0.5$, $L_I = 15$, $L_R = 150$. Solid orange curves are $[I](t)$ under the intervention, while dashed orange curves are $[I](t)$ without any intervention. Gray dashed lines denote the starts of the intervention, holding, and relaxation periods.}
    \label{fig:Examples1}
\end{figure}

\begin{figure}[htp]
    \centering
    \begin{subfigure}{\textwidth}
    \centering
    \includegraphics[width=1\textwidth]{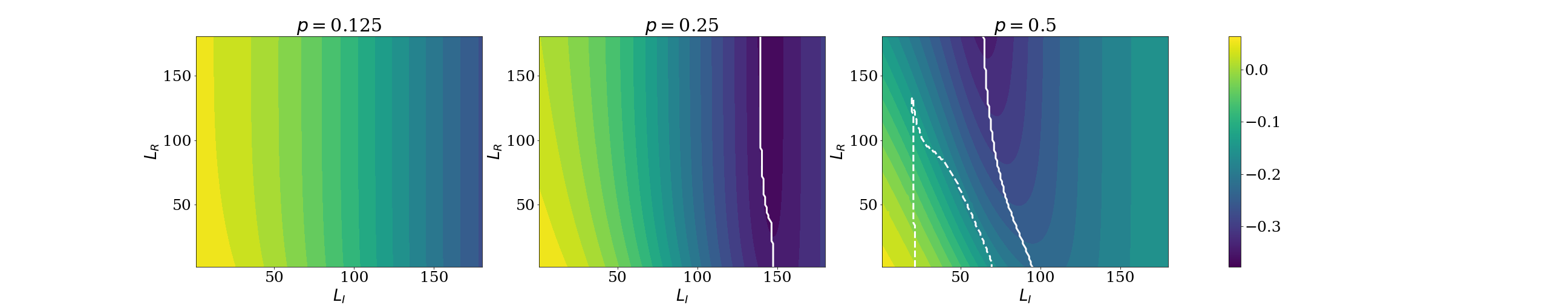}
    \caption{}
    \label{fig:int_1a}
    \end{subfigure}
    \hfill
    \begin{subfigure}{\textwidth}
    \centering
    \includegraphics[width=1\textwidth]{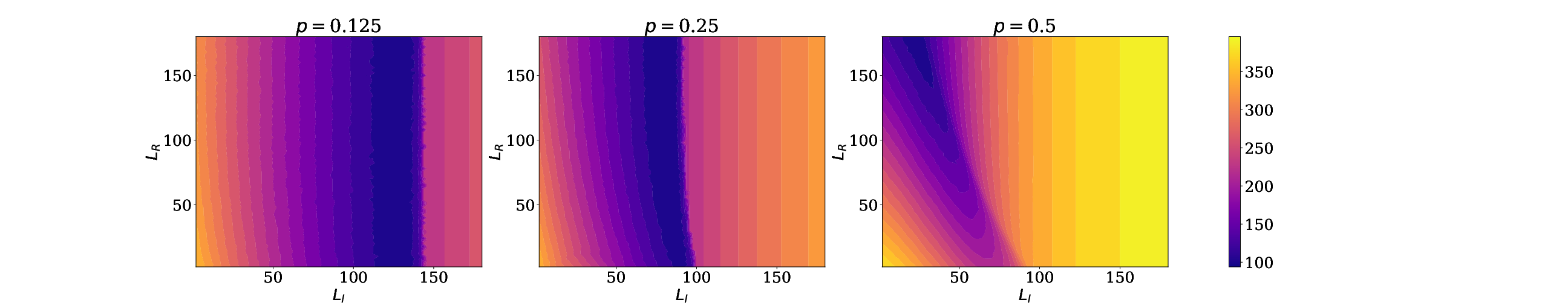}
    \caption{}
    \label{fig:int_1b}
    \end{subfigure}
    \caption{Plots of the RCFS (a) and AIAT (b) for the $L_H = 15$ and $q=0.01$. For intensities $p=0.125, 0.25,$ and $0.5$, the intervention period and relaxation period lengths $L_I$ and $L_R$ vary from $2$ to $180$ days. In (a), the solid white curve denotes the qualitative boundary, to the right of which uniform spikes occur. The dashed white line in the third panel denotes the boundary of the region where two spikes occur.}
    \label{fig:int_1}
\end{figure}

To achieve this effect in the evolution of the average number of contacts, we consider link activation rate $\alpha(t)$ and deletion rate $\omega(t)$ functions that are piecewise constant. These rate functions can be seen in the bottom two panels of Fig. \ref{fig:intervention1}. Since contacts are only broken in the first phase, $\omega(t) = \omega^*$ in the first phase and $0$ otherwise. Since contacts are only formed in the third phase, $\alpha(t) = \alpha^*$ in the third phase and $0$ otherwise. As the dynamical equation for $\langle k \rangle$ (\ref{eq:k}) is a first-order linear ODE, the resulting curve for $\langle k \rangle(t)$ will be piecewise exponential, and the values of $\alpha^*$ and $\omega^*$ are easily computed for a given $p.$ Other than $p,$ four other parameters characterize the simple intervention: the lengths of the three phases $L_I,L_H,$ and $L_R$, and the threshold proportion of the population $q\in [0,1)$ to initiate the intervention. The full simple intervention scheme can be described as follows:  
\begin{itemize}
    \item No intervention: the epidemic spreads unabated until $[I]$ increases through $qN$ ($\alpha =\omega = 0$).
    \item Intervention Phase (length $L_I$): intervention occurs, edges are removed at a constant rate ($\alpha = 0,\omega = \omega^*$).
    \item Holding Phase (length $L_H$): intervention holds, edges are neither removed nor added ($\alpha = \omega =0$).
    \item Relaxation Phase (length $L_R$): interventions are relaxed, edges are added at a constant rate ($\alpha = \alpha^*,\omega =0)$.
\end{itemize}
As this scheme requires five ``intervention" parameters, $p,q,L_I,L_H,$ and $L_R$, exploring the full impact of the interventions is difficult. To better see the effects, we consider an example scheme where we fix two parameter values in each and allow the other three to vary. To focus on the impact of the severity parameter $p$ and the lengths of the intervention and relaxation phases $L_I$ and $L_R$, we set $L_H = 15$ and $q = 0.01$ for the remainder of this section. Thus, the intervention begins when infections reach one percent of the population, and the holding phase is fixed at 15 days for all interventions. The other three parameters are allowed to vary. This allows for both abrupt and gradual implementations of interventions and relaxation of measures, and different levels of intervention intensity. Fig. \ref{fig:Examples1} shows the prevalence of some example intervention schemes, showing rich qualitative behavior.
To assess the effectiveness of the simple intervention we plot the RCFS and the AIAT for a large number of parameter combinations. We allow the lengths of both the intervention and relaxation periods $L_I$ and $L_R$ to vary from $2$ to $180$ days, and consider three different intensities $p = 0.125,0.25,0.5.$ The results are shown in Fig. \ref{fig:int_1}.

\begin{figure}[htp]
    \centering
    \begin{subfigure}[b]{0.3\textwidth}
    \centering
    \includegraphics[width=\textwidth]{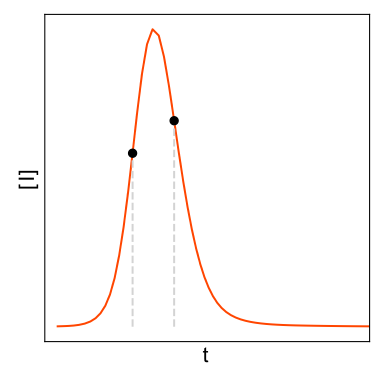}
    \caption{}
    \label{fig:typea}
    \end{subfigure}
    \begin{subfigure}[b]{0.3\textwidth}
    \centering
    \includegraphics[width=\textwidth]{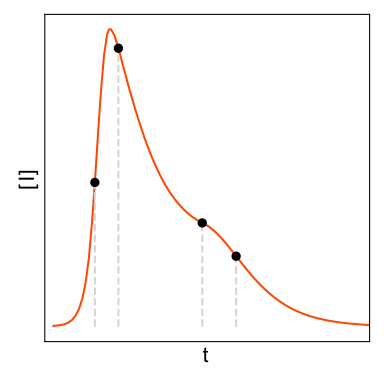}
    \caption{}
    \label{fig:typeb}
    \end{subfigure}
    \begin{subfigure}[b]{0.3\textwidth}
    \centering
    \includegraphics[width=\textwidth]{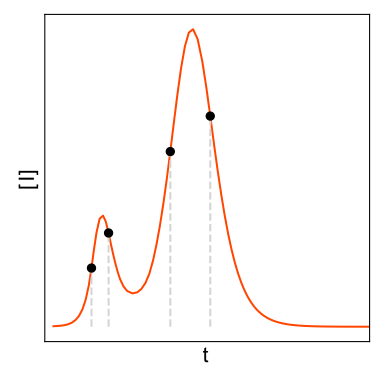}
    \caption{}
    \label{fig:typec}
    \end{subfigure}
    \caption{Types of infection curves with the simple intervention: (a) uniform spike, (b) non-uniform spike, (c) multiple spikes. Black dots denote inflection points.}
    \label{fig:type}
\end{figure}

A significant common feature of the plots in Fig. \ref{fig:int_1a} is a qualitative boundary (solid white) that divides $(L_I,L_R)$ space into two distinct classes of the resulting infection curve (for $p=0.125$, this occurs outside the boundaries of the plot). To the right of the boundary, infection curves are characterized by a single ``uniform spike," defined by an prevalence curve $[I](t)$ with two inflection points and a single local maximum (Fig. \ref{fig:typea}). To the left of the boundary, infection curves take the form of either a single ``non-uniform spike" (Fig. \ref{fig:typeb}), with more than two inflection points but only one local maximum, or multiple spikes (Fig. \ref{fig:typec}), with more than two inflection points and multiple local maxima. For $p=0.125$ and $p=0.25,$ only multiple spikes occur to the left of the boundary. For small $L_I$, the first spike is small and the second spike is large, and occasionally the final size of the epidemic surpasses the static case due to network alterations. As $L_I$ approaches the qualitative boundary, the second spike becomes shorter and occurs later until negligible. This phenomenon can also be seen in Fig. \ref{fig:int_1b}: as $L_I$ increases, the AIAT decreases until the second spike drops below the threshold $qN,$ at which point the AIAT increases as the first spike grows taller. For $p=0.5$ on the other hand, both nonuniform spikes and multiple spikes are possible to the left of the boundary. Multiple spikes occur in the region of $(L_I,L_R)$ space enclosed by the dashed white curve, while a single nonuniform spike occurs elsewhere left of the qualitative boundary. 

A few other observations warrant comment. First, the length of the intervention $L_I$ appears to be more important in determining epidemic's final size compared to $L_R$. This is intuitive, as the most significant changes to network structure occur during the intervention phase. Second, as $p$ increases, the qualitative boundary shifts generally left. This means that for less severe interventions, single uniform spikes will occur for smaller $L_I$ values. This observation carries weight for repeated interventions, explored in Section \ref{prevdep}, as single uniform spikes are heavily penalized by the AIAT. Third, nonuniform spikes occur for $p=0.5,$ but not for $p=0.125$ or $p=0.25.$ We hypothesize that there may exist some threshold $p^*$ where nonuniform spikes don't occur below $p^*,$ but do above $p^*.$

\subsection{Prevalence-Dependent Intervention} \label{prevdep}

While the simple intervention scheme provides a simple yet general model of social distancing, its implementation lacks a degree of realism. Interventions are put into place only once, and the epidemic continues, often with infections spiking after measures begin to relax. In reality, we would expect public health measures to be responsive to rising prevalence. Moreover, continued interventions might be triggered by some indicator, such as case numbers, deaths, hospital capacity, etc... In this section, we adapt the intervention scheme from Section \ref{simple} so that it may be reimplemented when a prevalence-based condition is satisfied, forming the prevalence-dependent intervention.
\begin{figure}[h]
    \centering
    \includegraphics[width=\textwidth]{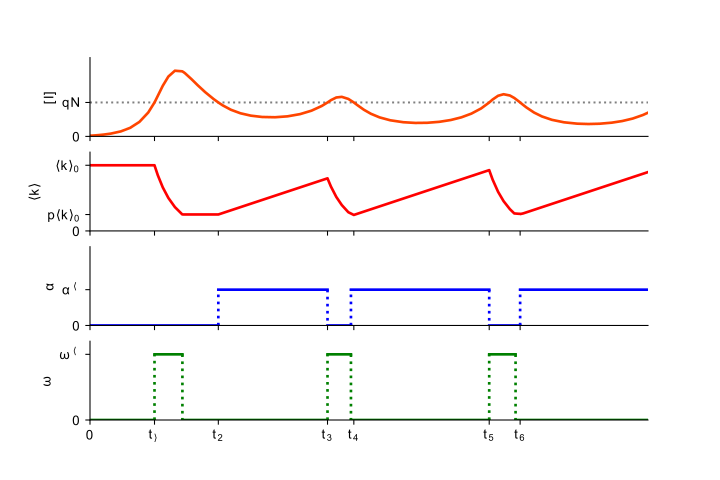}
    \caption{Prevalence-Dependent Intervention. The intervention begins when $[I] = qN,$ and edges are deleted at a constant rate $\omega^*$ until $\langle k \rangle$ decreases to $p\langle k \rangle_0,$ at which point there is no change to the network until $[I]$ drops below the threshold $qN.$ Then, edges are added at a constant rate $\alpha^*$ until $\langle k\rangle$ returns to $\langle k \rangle_0$ or $[I]$ increases through the threshold $qN$, at which point the intervention begins again.}
    \label{fig:prevalencedependentinterventionscheme}
\end{figure}
We begin with two more realistic assumptions about how a public health response might unfold. First, interventions are reimplemented any time the prevalence increases through some threshold. Second, the relaxation phase of an intervention doesn't begin until the prevalence has dropped below the threshold. We incorporate these assumptions into a new prevalence-dependent intervention scheme. The scheme is determined by four parameters: $q,p,L_I,$ and $L_R.$ As before, interventions begin when $[I]$ reaches $qN$, $p$ is the severity of the intervention, and $L_I$ and $L_R$ are now the maximum lengths of the intervention and relaxation periods, which determine $\omega^*$ and $\alpha^*$ as in Section \ref{simple}. We can define the new scheme as follows:
\begin{itemize}
    \item As $[I]$ increases through $qN,$ a new intervention is implemented.
    \item Intervention Phase: Once an intervention is implemented, edges are deleted at rate $\omega = \omega^*$ until $\langle k \rangle = p\langle k \rangle_0.$
    \item Holding Phase: At the end of the intervention period, a holding period begins ($\alpha=\omega = 0$) until the prevalence has dropped below the threshold $qN.$ If the prevalence drops below the threshold during the intervention period, the holding period has length 0.
    \item Relaxation Phase: Edges are added at rate $\alpha = \alpha^*$ until $\langle k \rangle = \langle k \rangle_0,$ or a new intervention is implemented.
\end{itemize}
It worth noting that compared to the simple intervention in Section \ref{simple}, the intervention, holding, and relaxation phases can all be of variable length. For instance, if the average number of contacts $\langle k \rangle$ has not rebounded to $\langle k \rangle_0$ by the time a new implementation begins, the resulting relaxation period is shorter than $L_R$. Moreover, in the subsequent intervention phase, edges delete until $\langle k \rangle = p\langle k \rangle_0$ and the phase is shorter than $L_I.$ In sum, while $\omega^*$ and $\alpha^*$ are fixed, the average number of contacts is never less than $p\langle k \rangle_0$ and the effective lengths of different intervention and relaxation phases may vary. An example implementation of the prevalence-dependent scheme is shown in Fig. \ref{fig:prevalencedependentinterventionscheme}, which shows both holding periods of nonzero length as well as intervention and relaxation periods that are shorter than $L_I$ and $L_R$ respectively.

\begin{figure}[h]
    \centering
    \begin{subfigure}[b]{.3\textwidth}
    \centering
    \includegraphics[width=\textwidth]{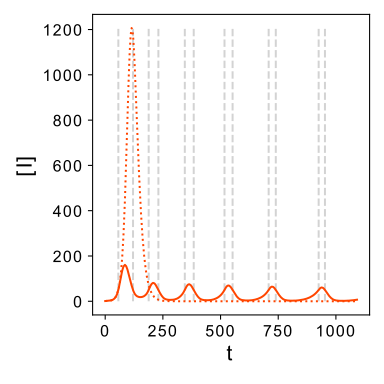}
    \caption{}
    \end{subfigure}
    \begin{subfigure}[b]{.3\textwidth}
    \centering
    \includegraphics[width=\textwidth]{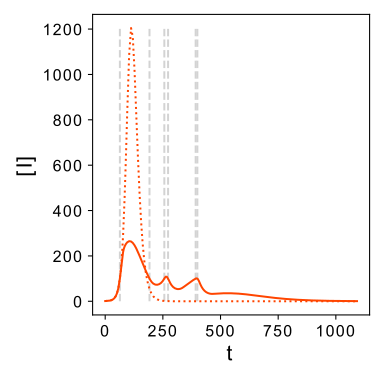}
    \caption{}
    \end{subfigure}
    \begin{subfigure}[b]{0.3\textwidth}
    \centering
    \includegraphics[width=\textwidth]{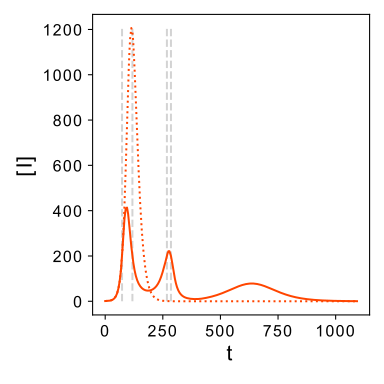}
    \caption{}
    \end{subfigure}
    \caption{Example infection curves $[I](t)$ for the prevalence-dependent intervention. Parameters shown are (a) $q= 0.005, p = 0.125, L_I = 60, L_R = 60$, (b) $q=0.01, p = 0.5, L_I = 15, L_R = 60$, (c) $q=0.02, p = 0.25, L_I = 30, L_R = 120$. Solid orange curves are $[I](t)$ under the intervention, while dashed orange curves are $[I](t)$ without any intervention. Dashed gray lines denote times when $[I] = qN$.}
    \label{fig:PDexample}
\end{figure}
\begin{figure}[h]
    \centering
    \includegraphics[width=\textwidth]{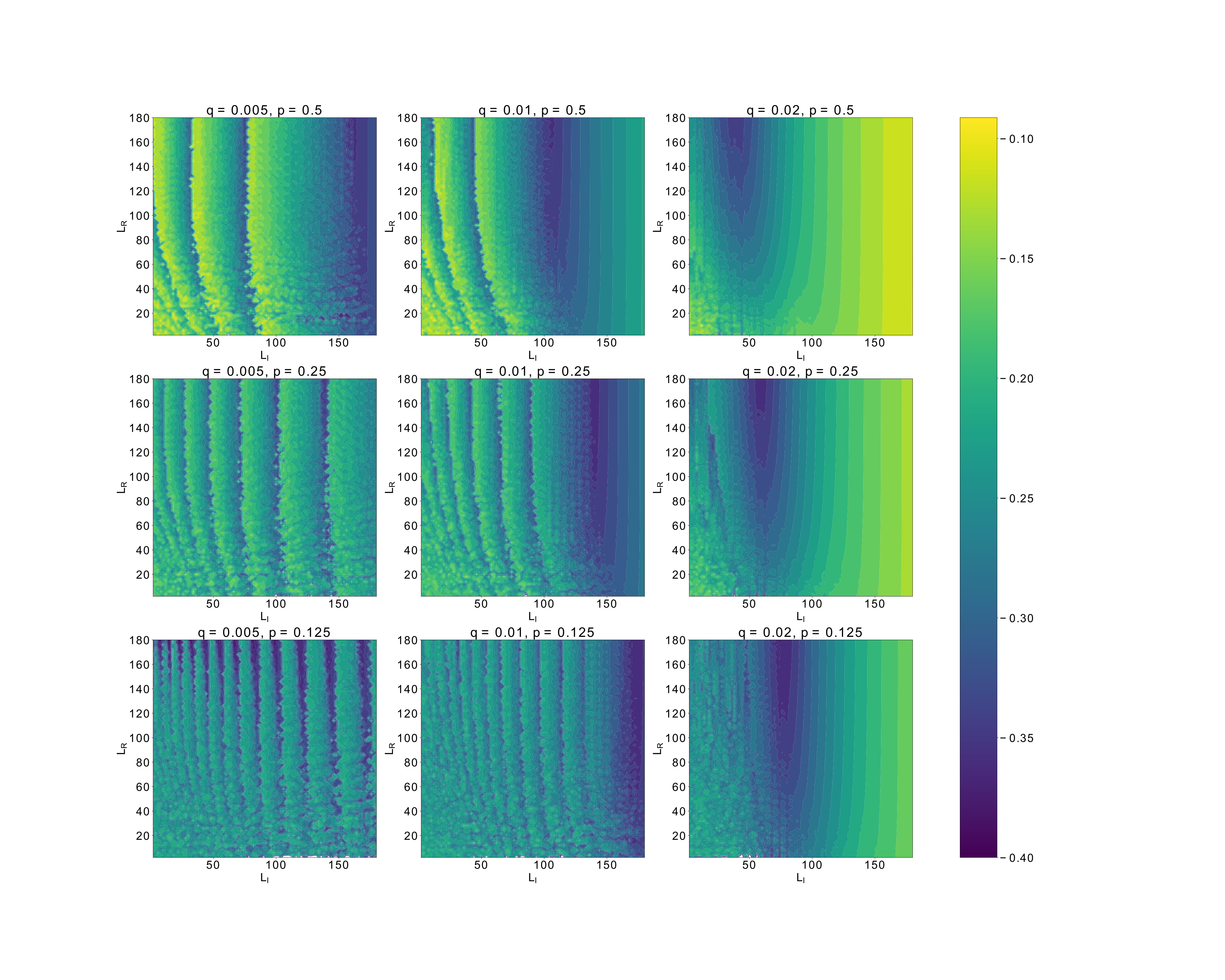}
    \caption{Relative change in final size (RCFS) for the prevalence-dependent intervention. Each plot represents a choice of $p$ and $q,$ with $L_I$ and $L_R$ on the axes, ranging from $2$ for $180$.}
    \label{fig:PDRCFS}
\end{figure}
A notable feature of the prevalence-dependent intervention is its ability to generate infection curves with multiple spikes as the epidemic progresses. Examples of this behavior are shown in Fig. \ref{fig:PDexample}. To fully explore the intervention, we again consider the RCFS for a variety of parameter combinations. Fig. \ref{fig:PDRCFS} shows the RCFS for different thresholds ($q = 0.005,0.01,0.02$) and intensities ($p = 0.125,0.25,0.5$) as $L_I$ and $L_R$ both vary from $2$ to $180$ days. Though not shown, as with the simple intervention each case has a qualitative boundary, to the right of which infection curves are single, uniform spikes. The most significant departure from the simple intervention though is to the left of the qualitative boundary. In the simple case, infection curves from this region took the form of either two spikes or a single nonuniform spike. With the prevalence-dependent intervention, the infection curve behavior is richer. 

\begin{figure}[h]
    \centering
    \begin{subfigure}[b]{.3\textwidth}
    \centering
    \includegraphics[width=\textwidth]{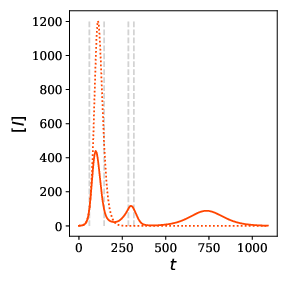}
    \caption{}
    \label{fig:ExamplePDa}
    \end{subfigure}
    \begin{subfigure}[b]{.3\textwidth}
    \centering
    \includegraphics[width=\textwidth]{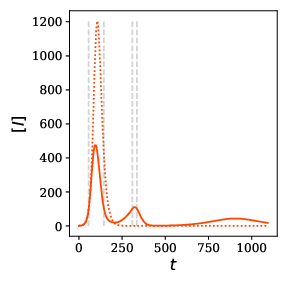}
    \caption{}
    \label{fig:ExamplePDb}
    \end{subfigure}
    \begin{subfigure}[b]{0.3\textwidth}
    \centering
    \includegraphics[width=\textwidth]{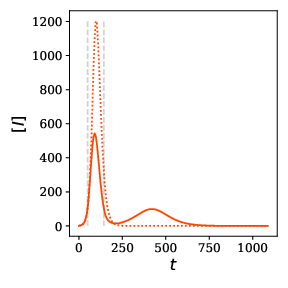}
    \caption{}
    \label{fig:ExamplePDc}
    \end{subfigure}
    \caption{Progression of the infection curve $[I](t)$ as $L_I$ increases, showing the shrinking of the final spike and the penultimate spike dropping below the threshold $qN$. Parameters are $q=0.01,p=0.25,L_R = 90$ and $L_I = 70$ (a), $78$ (b), $92$ (c). Solid orange curves are $[I](t)$ under the intervention, while dashed orange curves are $[I](t)$ without any intervention. Dashed gray lines denote times when $[I] = qN$.}
    \label{fig:ExamplesPD}
\end{figure}

\begin{figure}[h]
    \centering
    \includegraphics[width=\textwidth]{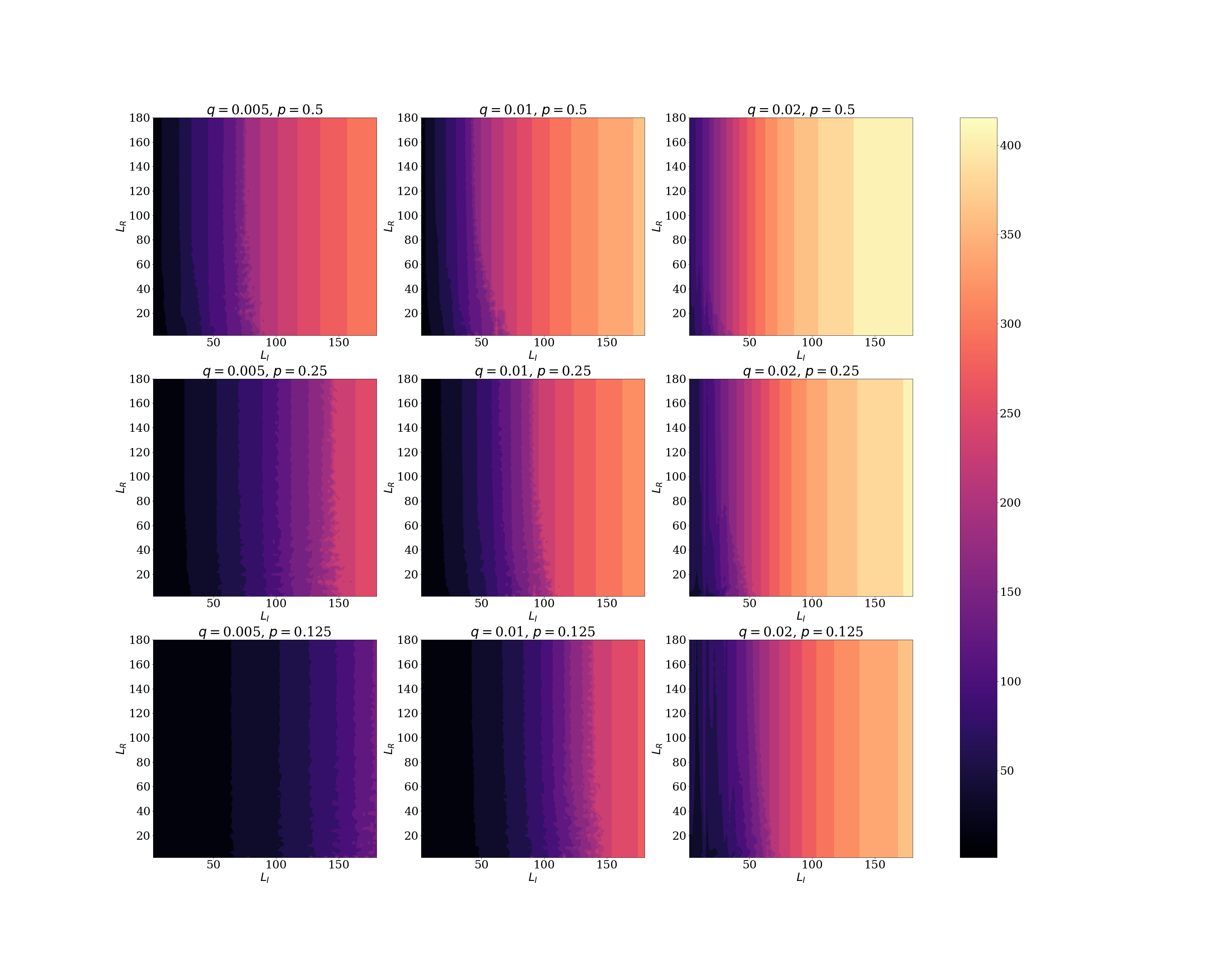}
    \caption{Average infections above threshold (AIAT) for the prevalence-dependent intervention. Each plot represents a choice of $p$ and $q,$ with $L_I$ and $L_R$ on the axes, ranging from $2$ for $180$.}
    \label{fig:PDAIAT}
\end{figure}

The region is characterized by ``waves" in the RCFS, particularly for lower values of $p$. The boundaries of these waves can be described by the number of spikes that occur over the course of the epidemic. Holding $L_R$ fixed and increasing $L_I$ through one of these contours helps explain the behavior of the infection curve in this region (Fig. \ref{fig:ExamplesPD}). At the crest, the final spike peaks just below the threshold $qN$ (Fig. \ref{fig:ExamplePDa}). As $L_I$ increases, the final spike occurs later and peaks lower (Fig. \ref{fig:ExamplePDb}) and the RCFS decreases until the spike vanishes. Then, the penultimate spike becomes the new final spike, peaking just below the threshold (Fig. \ref{fig:ExamplePDc}) and the RCFS jumps up as a new wave crests. This underscores a potential limitation of a threshold-based intervention: if a spike does not reach the threshold and no intervention occurs, the spike occurs over a longer period of time and more infections accumulate than if the spike had triggered an intervention. A practical implication of this observation is that no spike in infections should go unaddressed by interventions if the goal is to reduce the number of cumulative infections. We also consider the AIAT for the same parameter combinations (Fig. \ref{fig:PDAIAT}), though the conclusions by this metric are less complex. For any combination of $p$ and $q,$ increasing $L_I$ leads to a larger AIAT. This suggests that when considering interventions with the same RCFS, more abrupt interventions (smaller $L_I$) are preferable. However, an interesting observation is that the AIAT increases rapidly as the epidemic changes from three to two spikes.
\begin{figure}[h]
    \centering
    \begin{subfigure}[b]{\textwidth}
    \centering
    \includegraphics[width=\textwidth]{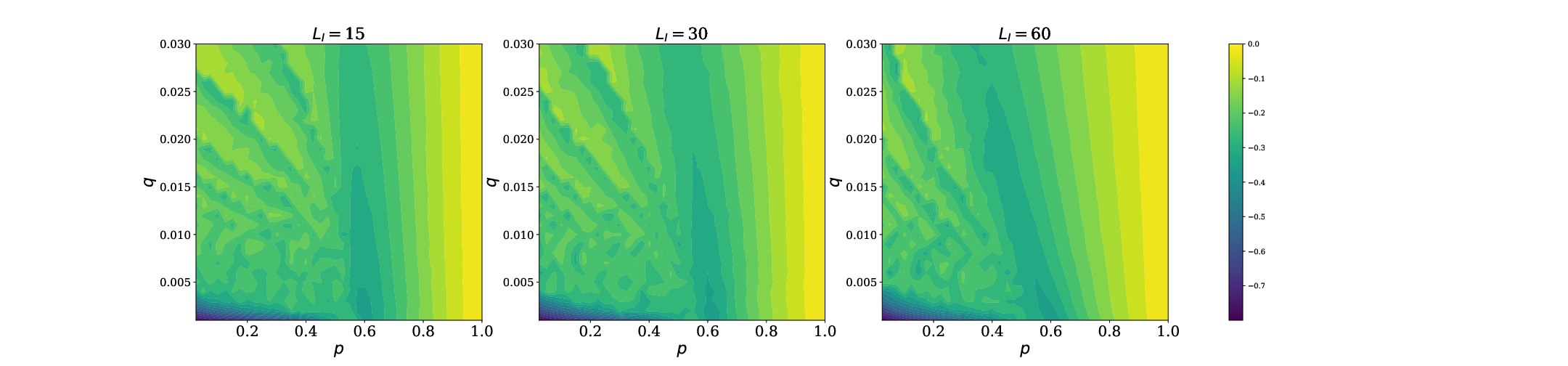}
    \caption{}
    \end{subfigure}
    \begin{subfigure}[b]{\textwidth}
    \centering
    \includegraphics[width=\textwidth]{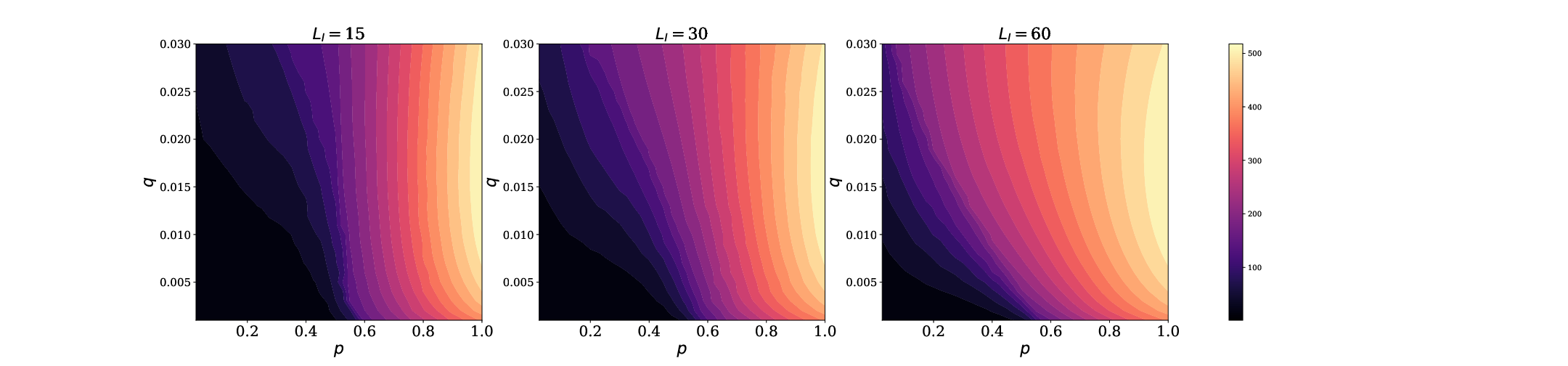}
    \caption{}
    \end{subfigure}
    \caption{Plots of the RCFS (a) and AIAT (b) for the prevalence-dependent intervention with $L_I = 15,30,60$ and $L_R = 90$ as $p$ varies from 0 to 1 and $q$ varies from $0$ to $0.03$. Notably, both measures indicate highly-effective interventions for small values of $p$ and $q$.}
    \label{fig:PQPlots}
\end{figure}

While Figs. \ref{fig:PDRCFS} and \ref{fig:PDAIAT} show the overall behavior of the prevalence-dependent intervention, by considering fixed values of $L_I$ and $L_R$ and allowing $p$ and $q$ to vary, we get a more pointed perspective on the effectiveness of this type of intervention. Fig. \ref{fig:PQPlots} shows increasingly gradual interventions from left to right with plots of the RCFS and AIAT as $p$ and $q$ vary on the axes. Notably, regardless of $L_I,$ low values of $p$ and $q$ are able to produce interventions that both greatly decrease the final size of the epidemic, and the average infections above threshold. This suggests that for sufficiently low thresholds ($q$) and sufficiently severe decreases in contacts ($p$), the length over which the decrease in contacts occurs ($L_I$) does not play an important role in the effectiveness of interventions. However, as $q$ or $p$ increases, $L_I$ has a more pronounced impact. In particular, for low values of $p$ and large values of $q,$ a longer, more gradual intervention can lead to more average infections above threshold. Moreover, a stark change in both effectiveness metrics occurs for large values of $p,$ (around $p=0.5$ for $L_I=15$ and $L_I=30$). This suggests that if an intervention doesn't reduce average contacts sufficiently, a highly effective intervention isn't possible, regardless of the other parameter values.

\section{Discussion} \label{discussion}

In this paper, we have developed a new SEIR model on a network with random link activation/deletion dynamics. Using piecewise constant activation and deletion rate functions, we propose two simple mechanisms for social distancing interventions. The simple intervention models a single intervention event, where contacts are decreased over a period of time, stay constant, and then return to pre-intervention levels. The prevalence-dependent intervention expands the simple case to more complex scenarios, where interventions can be reintroduced in the face of rising prevalence. Using the unipartite projection of a bipartite network, and epidemiological parameters representative of COVID-19, we examine the effectiveness of a wide range of potential social distancing policies on relatively large heterogeneous, clustered networks.

Both intervention schemes are shown to capture a wide variety of behaviors in the prevalence ``curve," which has received considerable attention in both academic studies and public health messaging. The simple intervention manifests curves with one or two spikes, while the curves for prevalence-dependent intervention can have many more. Moreover, the behavior of the prevalence curve is consistent across a number of parameters and can be described qualitatively with success. This is despite the simplicity of social distancing mechanism introduced by the piecewise constant activation and deletion rates $\alpha(t)$ and $\omega(t),$ which take on values $\alpha^*$ or $\omega^*$ respectively, or zero. We have not considered the cases where the values of $\alpha^*$ and $\omega^*$ may change over time, or where $\alpha(t)$ and $\omega(t)$ are not piecewise constant. As such, our model has natural extensions that may capture an even richer variety of qualitative behaviors. 

Furthermore, the mechanisms proposed in this paper offer insights into what makes for a successful intervention. We have used two metrics as simplified public health goals to evaluate the effectiveness of interventions: the relative change in final size (RCFS) and the average infections above threshold (AIAT). For the more realistic prevalence-dependent intervention scheme, we find that the most effective interventions come when the threshold number of infections is low and the intervention severely decreases average contacts. When these conditions are met, the relative change in the final size is greatly decreased and the length over which the intervention is implemented has little impact on the effectiveness. However, even small increases in the threshold value can greatly impact the effectiveness of interventions. As well, if interventions do not sufficiently reduce contacts (around fifty percent), they are rendered significantly less effective by both measures.

While this is a first foray into the use of adaptive networks to model social distancing for an SEIR disease, we acknowledge some limitations of our model. First, there is a trade-off between complexity of the disease natural history model and the number of equations of the pairwise model; age-structured models or other more complex compartmental models are popular for COVID-19, but added compartments require tracking an increasing number of edge types. However, even simple extensions (such as the inclusion of an asymptomatic infectious state) present interesting opportunities. Second, while the random link activation/deletion process is simple to implement, it has some unrealistic features. In particular, in the $t\to \infty$ limit, one can show from the degree distribution generating function that the resulting network approaches an Erd\H{o}s-R\'{e}nyi random graph, with vanishing clustering and an approximately Poisson degree distribution. One manifestation of this property is a rapidly declining clustering coefficient over time. While the piecewise constant activation and deletion rates mitigate this to an extent, the network resulting from these social distancing policies is fundamentally different than the initial network state. To overcome this limitation, future investigations might involve new processes for network dynamics, such as activation/deletion on a fixed network or network dynamics on an underlying bipartite mixing network. 
\section*{Acknowledgements}
This paper is a continuation of a project that began at the ``Dynamics and data in the COVID-19 pandemic" workshop hosted by the American Institute of Mathematics. The authors would like to thank Stephen Schecter, Hans Kaper, and the rest of the workshop staff for their guidance. We also thank Alan Hastings for his insightful comments.

\bibliographystyle{spbasic}      
\bibliography{references.bib}   

\begin{thebibliography}{31}
\providecommand{\natexlab}[1]{#1}
\providecommand{\url}[1]{{#1}}
\providecommand{\urlprefix}{URL }
\expandafter\ifx\csname urlstyle\endcsname\relax
  \providecommand{\doi}[1]{DOI~\discretionary{}{}{}#1}\else
  \providecommand{\doi}{DOI~\discretionary{}{}{}\begingroup
  \urlstyle{rm}\Url}\fi
\providecommand{\eprint}[2][]{\url{#2}}

\bibitem[{Ahmed et~al(2018)Ahmed, Zviedrite, and
  Uzicanin}]{ahmed_effectiveness_2018}
Ahmed F, Zviedrite N, Uzicanin A (2018) Effectiveness of workplace social
  distancing measures in reducing influenza transmission: A systematic review.
  BMC Public Health 18(1):518

\bibitem[{Anastassopoulou et~al(2020)Anastassopoulou, Russo, Tsakris, and
  Siettos}]{anastassopoulou_data-based_2020}
Anastassopoulou C, Russo L, Tsakris A, Siettos C (2020) Data-based analysis,
  modelling and forecasting of the {{COVID}}-19 outbreak. PLOS ONE
  15(3):e0230,405

\bibitem[{Chang et~al(2021)Chang, Pierson, Koh, Gerardin, Redbird, Grusky, and
  Leskovec}]{chang_mobility_2021}
Chang S, Pierson E, Koh PW, Gerardin J, Redbird B, Grusky D, Leskovec J (2021)
  Mobility network models of {{COVID}}-19 explain inequities and inform
  reopening. Nature 589(7840):82--87

\bibitem[{Davey et~al(2008)Davey, Glass, Min, Beyeler, and
  Glass}]{davey_effective_2008}
Davey VJ, Glass RJ, Min HJ, Beyeler WE, Glass LM (2008) Effective, {{Robust
  Design}} of {{Community Mitigation}} for {{Pandemic Influenza}}: {{A
  Systematic Examination}} of {{Proposed US Guidance}}. PLoS ONE 3(7):e2606

\bibitem[{Eames and Keeling(2002)}]{eames_modeling_2002}
Eames KTD, Keeling MJ (2002) Modeling dynamic and network heterogeneities in
  the spread of sexually transmitted diseases. Proceedings of the National
  Academy of Sciences 99(20):13,330--13,335

\bibitem[{Eubank et~al(2004)Eubank, Guclu, Anil~Kumar, Marathe, Srinivasan,
  Toroczkai, and Wang}]{eubank_modelling_2004}
Eubank S, Guclu H, Anil~Kumar VS, Marathe MV, Srinivasan A, Toroczkai Z, Wang N
  (2004) Modelling disease outbreaks in realistic urban social networks. Nature
  429(6988):180--184

\bibitem[{Eubank et~al(2020)Eubank, Eckstrand, Lewis, Venkatramanan, Marathe,
  and Barrett}]{eubank_commentary_2020}
Eubank S, Eckstrand I, Lewis B, Venkatramanan S, Marathe M, Barrett CL (2020)
  Commentary on {{Ferguson}}, et al., ``{{Impact}} of {{Non}}-pharmaceutical
  {{Interventions}} ({{NPIs}}) to {{Reduce COVID}}-19 {{Mortality}} and
  {{Healthcare Demand}}''. Bulletin of Mathematical Biology 82(4):52

\bibitem[{Ferguson et~al(2020)Ferguson, Laydon, Nedjati~Gilani, Imai, Ainslie,
  Baguelin, Bhatia, Boonyasiri, Cucunuba~Perez, {Cuomo-Dannenburg}, Dighe,
  Dorigatti, Fu, Gaythorpe, Green, Hamlet, Hinsley, Okell, Van~Elsland,
  Thompson, Verity, Volz, Wang, Wang, Walker, Winskill, Whittaker, Donnelly,
  Riley, and Ghani}]{ferguson_report_2020}
Ferguson N, Laydon D, Nedjati~Gilani G, Imai N, Ainslie K, Baguelin M, Bhatia
  S, Boonyasiri A, Cucunuba~Perez Z, {Cuomo-Dannenburg} G, Dighe A, Dorigatti
  I, Fu H, Gaythorpe K, Green W, Hamlet A, Hinsley W, Okell L, Van~Elsland S,
  Thompson H, Verity R, Volz E, Wang H, Wang Y, Walker P, Winskill P, Whittaker
  C, Donnelly C, Riley S, Ghani A (2020) Report 9: {{Impact}} of
  non-pharmaceutical interventions ({{NPIs}}) to reduce {{COVID19}} mortality
  and healthcare demand. Tech. rep., {Imperial College London}

\bibitem[{Glass et~al(2006)Glass, Glass, Beyeler, and
  Min}]{glass_targeted_2006}
Glass RJ, Glass LM, Beyeler WE, Min HJ (2006) Targeted {{Social Distancing
  Design}} for {{Pandemic Influenza}}. Emerging Infectious Diseases 12(11):11

\bibitem[{Gross and Sayama(2009)}]{gross_adaptive_2009}
Gross T, Sayama H (eds)  (2009) Adaptive {{Networks}}: {{Theory}}, {{Models}}
  and {{Applications}}. Understanding {{Complex Systems}}, {Springer Berlin
  Heidelberg}, {Berlin, Heidelberg}

\bibitem[{Gross et~al(2006)Gross, D'Lima, and Blasius}]{gross_epidemic_2006}
Gross T, D'Lima CJD, Blasius B (2006) Epidemic {{Dynamics}} on an {{Adaptive
  Network}}. Physical Review Letters 96(20):208,701

\bibitem[{House and Keeling(2011)}]{house_insights_2011}
House T, Keeling MJ (2011) Insights from unifying modern approximations to
  infections on networks. Journal of The Royal Society Interface 8(54):67--73

\bibitem[{Keeling(1999)}]{keeling_effects_1999-1}
Keeling MJ (1999) The effects of local spatial structure on epidemiological
  invasions. Proceedings Biological sciences 266(1421):859--867

\bibitem[{Keeling et~al(1997)Keeling, Rand, and
  Morris}]{keeling_correlation_1997}
Keeling MJ, Rand DA, Morris AJ (1997) Correlation models for childhood
  epidemics. Proceedings of the Royal Society of London Series B: Biological
  Sciences 264(1385):1149--1156

\bibitem[{Kiss et~al(2012)Kiss, Berthouze, Taylor, and
  Simon}]{kiss_modelling_2012}
Kiss IZ, Berthouze L, Taylor TJ, Simon PL (2012) Modelling approaches for
  simple dynamic networks and applications to disease transmission models.
  Proceedings of the Royal Society A: Mathematical, Physical and Engineering
  Sciences 468(2141):1332--1355

\bibitem[{Kiss et~al(2017)Kiss, Miller, and Simon}]{kiss_mathematics_2017}
Kiss IZ, Miller JC, Simon PL (2017) Mathematics of {{Epidemics}} on
  {{Networks}}: {{From Exact}} to {{Approximate Models}}, Interdisciplinary
  {{Applied Mathematics}}, vol~46. {Springer International Publishing}, {Cham}

\bibitem[{Li et~al(2020)Li, Guan, Wu, Wang, Zhou, Tong, Ren, Leung, Lau, Wong,
  Xing, Xiang, Wu, Li, Chen, Li, Liu, Zhao, Liu, Tu, Chen, Jin, Yang, Wang,
  Zhou, Wang, Liu, Luo, Liu, Shao, Li, Tao, Yang, Deng, Liu, Ma, Zhang, Shi,
  Lam, Wu, Gao, Cowling, Yang, Leung, and Feng}]{li_early_2020}
Li Q, Guan X, Wu P, Wang X, Zhou L, Tong Y, Ren R, Leung KS, Lau EH, Wong JY,
  Xing X, Xiang N, Wu Y, Li C, Chen Q, Li D, Liu T, Zhao J, Liu M, Tu W, Chen
  C, Jin L, Yang R, Wang Q, Zhou S, Wang R, Liu H, Luo Y, Liu Y, Shao G, Li H,
  Tao Z, Yang Y, Deng Z, Liu B, Ma Z, Zhang Y, Shi G, Lam TT, Wu JT, Gao GF,
  Cowling BJ, Yang B, Leung GM, Feng Z (2020) Early {{Transmission Dynamics}}
  in {{Wuhan}}, {{China}}, of {{Novel Coronavirus}}\textendash{{Infected
  Pneumonia}}. New England Journal of Medicine 382(13):1199--1207

\bibitem[{Linton et~al(2020)Linton, Kobayashi, Yang, Hayashi, Akhmetzhanov,
  Jung, Yuan, Kinoshita, and Nishiura}]{linton_incubation_2020}
Linton N, Kobayashi T, Yang Y, Hayashi K, Akhmetzhanov A, Jung Sm, Yuan B,
  Kinoshita R, Nishiura H (2020) Incubation {{Period}} and {{Other
  Epidemiological Characteristics}} of 2019 {{Novel Coronavirus Infections}}
  with {{Right Truncation}}: {{A Statistical Analysis}} of {{Publicly Available
  Case Data}}. Journal of Clinical Medicine 9(2):538

\bibitem[{Miller(2009)}]{miller_spread_2009}
Miller JC (2009) Spread of infectious disease through clustered populations.
  Journal of The Royal Society Interface 6(41):1121--1134

\bibitem[{Newman et~al(2001)Newman, Strogatz, and Watts}]{newman_random_2001}
Newman MEJ, Strogatz SH, Watts DJ (2001) Random graphs with arbitrary degree
  distributions and their applications. Physical Review E 64(2):026,118

\bibitem[{{Pastor-Satorras} et~al(2015){Pastor-Satorras}, Castellano,
  Van~Mieghem, and Vespignani}]{pastor-satorras_epidemic_2015}
{Pastor-Satorras} R, Castellano C, Van~Mieghem P, Vespignani A (2015) Epidemic
  processes in complex networks. Reviews of Modern Physics 87(3):925--979

\bibitem[{Rand(1999)}]{mcglade_correlation_1999}
Rand DA (1999) Correlation {{Equations}} and {{Pair Approximations}} for
  {{Spatial Ecologies}}. In: McGlade J (ed) Advanced {{Ecological Theory}},
  {Blackwell Publishing Ltd.}, {Oxford, UK}, pp 100--142

\bibitem[{Read et~al(2008)Read, Eames, and Edmunds}]{read_dynamic_2008}
Read JM, Eames KT, Edmunds WJ (2008) Dynamic social networks and the
  implications for the spread of infectious disease. Journal of The Royal
  Society Interface 5(26):1001--1007

\bibitem[{S{\'e}lley et~al(2015)S{\'e}lley, Besenyei, Kiss, and
  Simon}]{selley_dynamic_2015}
S{\'e}lley F, Besenyei {\'A}, Kiss IZ, Simon PL (2015) Dynamic {{Control}} of
  {{Modern}}, {{Network}}-{{Based Epidemic Models}}. SIAM Journal on Applied
  Dynamical Systems 14(1):168--187

\bibitem[{Shkarayev et~al(2014)Shkarayev, Tunc, and
  Shaw}]{shkarayev_epidemics_2014}
Shkarayev MS, Tunc I, Shaw LB (2014) Epidemics with temporary link deactivation
  in scale-free networks. Journal of Physics A: Mathematical and Theoretical
  47(45):455,006

\bibitem[{Taylor et~al(2012)Taylor, Simon, Green, House, and
  Kiss}]{taylor_markovian_2012}
Taylor M, Simon PL, Green DM, House T, Kiss IZ (2012) From {{Markovian}} to
  pairwise epidemic models and the performance of moment closure
  approximations. Journal of Mathematical Biology 64(6):1021--1042

\bibitem[{Tunc et~al(2013)Tunc, Shkarayev, and Shaw}]{tunc_epidemics_2013}
Tunc I, Shkarayev MS, Shaw LB (2013) Epidemics in {{Adaptive Social Networks}}
  with {{Temporary Link Deactivation}}. Journal of Statistical Physics
  151(1-2):355--366

\bibitem[{Valdez et~al(2012)Valdez, Macri, and
  Braunstein}]{valdez_intermittent_2012}
Valdez LD, Macri PA, Braunstein LA (2012) Intermittent social distancing
  strategy for epidemic control. Physical Review E 85(3):036,108

\bibitem[{You et~al(2020)You, Deng, Hu, Sun, Lin, Zhou, Pang, Zhang, Chen, and
  Zhou}]{you_estimation_2020}
You C, Deng Y, Hu W, Sun J, Lin Q, Zhou F, Pang CH, Zhang Y, Chen Z, Zhou XH
  (2020) Estimation of the time-varying reproduction number of {{COVID}}-19
  outbreak in {{China}}. International Journal of Hygiene and Environmental
  Health 228:113,555

\bibitem[{Youssef and Scoglio(2013)}]{youssef_mitigation_2013}
Youssef M, Scoglio C (2013) Mitigation of epidemics in contact networks through
  optimal contact adaptation. Mathematical Biosciences and Engineering
  10(4):1227--1251

\bibitem[{Zhang et~al(2020)Zhang, Litvinova, Wang, Wang, Deng, Chen, Li, Zheng,
  Yi, Chen, Wu, Liang, Wang, Yang, Sun, Longini, Halloran, Wu, Cowling, Merler,
  Viboud, Vespignani, Ajelli, and Yu}]{zhang_evolving_2020}
Zhang J, Litvinova M, Wang W, Wang Y, Deng X, Chen X, Li M, Zheng W, Yi L, Chen
  X, Wu Q, Liang Y, Wang X, Yang J, Sun K, Longini IM, Halloran ME, Wu P,
  Cowling BJ, Merler S, Viboud C, Vespignani A, Ajelli M, Yu H (2020) Evolving
  epidemiology and transmission dynamics of coronavirus disease 2019 outside
  {{Hubei}} province, {{China}}: A descriptive and modelling study. The Lancet
  Infectious Diseases 20(7):793--802

\end{thebibliography}

\newpage

\begin{appendices}

\numberwithin{equation}{section}
\numberwithin{figure}{section}

\section{Adaptive SEIR with Complex Closure}\label{AppA}

In this appendix, we develop an adaptive network SEIR pairwise model for heterogeneous, clustered networks. The model is analogous to the SIR model for heterogeneous, clustered networks  in \citet{house_insights_2011} with random link activation/deletion dynamics included.

The generic triple closure (\ref{eq:triple_approx}) proposed by \citet{house_insights_2011} can be further developed by introducing a new variable $\theta(t),$ the proportion of edges that have not transmitted the infection. With the observation that $[S_k] = Np_k\theta^k,$ (\ref{eq:triple_approx}) becomes
\begin{equation}
[ASI] \approx [AS][SI]\frac{g''(\theta)}{N(g'(\theta))^2}\left((1-\phi) + \phi g'(1) N \frac{[AI]}{\left(\sum_k k[A_k]\right) \left(\sum_k k [I_k]\right)}\right)\label{eq:Acomplexclosure}
\end{equation}
Moreover, we can express $\sum_k k[S_k] = N\theta g'(\theta),$ and we introduce auxiliary variables $Y = \sum_k k[E_k], Z = \sum_k k[I_k],$ and $\theta$. Observing that $\sum_k k [A_k] = [AS]+[AE]+[AI]+[AR],$ it follows that the dynamical equations for $Y$ and $Z$ (without network dynamics) are 
\begin{align}
    \dot{Y} &= \beta\frac{\theta g''(\theta)}{g'(\theta)}[SI] - \eta Y \label{eq:AY}\\
    \dot{Z} &= \eta Y - \gamma Z \label{eq:AZ}
\end{align}
Now we incorporate the effects of link activation and deletion. Notably, the probability generating function for the degree distribution is now time- dependent, taking the form 
\begin{equation}
    g(x,t) = \sum_{k=0}^{N-1}p_k(t)x^k.
\end{equation}
As a consequence, the ordinary derivatives of $g$ in (\ref{eq:Acomplexclosure})-(\ref{eq:AZ}) become partial derivatives with respect to $x.$ From (\ref{eq:[AA]}) and (\ref{eq:[AB]}), we can derive network- dynamical versions of (\ref{eq:AY}) and (\ref{eq:AZ}):
\begin{align}
    \dot{Y} &= \beta\frac{\theta g_{xx}(\theta,t)}{g_x(\theta,t)}[SI] - (\eta+\alpha+\omega)Y +\alpha(N-1)[E]\label{eq:AY2} \\
    \dot{Z} &= \eta Y - (\gamma+\alpha+\omega)Z+\alpha(N-1)[I] \label{eq:AZ2}
\end{align}
Next, the non-epidemiological network quantities in this model are entirely determined by the degree distribution probability generating function $g(x,t)$ and the clustering coefficient $\phi(t).$ We can express (\ref{eq:phi}) in terms of the generating function $g(x,t)$ as
\begin{equation}\label{eq:Aphi}
    \dot{\phi} = 3\alpha - \left(\alpha + \omega+2\alpha(N-2)\frac{g_x(1,t)}{g_{xx}(1,t)}\right)\phi.
\end{equation}
Equations (\ref{eq:Acomplexclosure}),(\ref{eq:AY2}), and (\ref{eq:Aphi}) require $g(x,t)$ and its derivatives explicitly, which can be found by solving (\ref{eq:pde}) using the method of characteristics:
\begin{equation}
    g(x,t) = g_0\left(\frac{\omega+\alpha x +\omega(x-1)e^{-(\alpha+\omega)t}}{\omega+\alpha x -\alpha(x-1)e^{-(\alpha+\omega)t}}\right)\left(\frac{\omega+\alpha x -\alpha(x-1)e^{-(\alpha+\omega)t}}{\alpha+\omega}\right)^{N-1}
\end{equation}
where $g_0(x) = g(x,0).$ Finally, we derive the evolution equation for $\theta(t)$ by differentiating $\dot{[S]} = Ng(\theta(t),t)$ and solving for $\dot\theta:$
\begin{equation}
    \dot{\theta} = -\frac{\beta[SI]}{Ng_x(\theta,t)}-(1-\theta)\left(\alpha\theta+\omega -\alpha(N-1)\frac{g(\theta,t)}{g_x(\theta,t)}\right).
\end{equation}
Thus, we arrive at the pairwise SEIR for a heterogeneous, clustered network with random link activation and deletion:
\begin{align}
[S] &= Ng(\theta,t)\\
\dot{[E]} &= \beta[SI]-\eta[E],\\
\dot{[I]} &= \eta[E]-\gamma[I],\\
\dot{[SS]} &= -2\beta[SSI]+\alpha[S]([S]-1)-(\alpha+\omega)[SS], \\
\dot{[SE]} &= \beta[SSI]-\beta[ESI] - \eta[SE]+\alpha[S][E]-(\alpha+\omega)[SE],\\
\dot{[SI]} &= \eta[SE] -\beta[SI] - \beta[ISI]-\gamma[SI]+\alpha[S][I]-(\alpha+\omega)[SI],\\
\dot{[EE]} &= 2\beta[ESI]-2\eta[EE]+\alpha[E]([E]-1)-(\alpha+\omega)[EE],\\
\dot{[EI]} &= \beta[ISI]+\beta[SI]+\eta[EE]-(\gamma+\eta)[EI]+\alpha[E][S]-(\alpha+\omega)[EI],\\
\dot{[II]} &= 2\eta[EI]-2\gamma[II]+\alpha[I]([I]-1)-(\alpha+\omega)[II],\\
\dot{Y} &= \beta\frac{\theta g_{xx}(\theta,t)}{g_x(\theta,t)}[SI] - (\eta+\alpha+\omega)Y +\alpha(N-1)[E],\\
\dot{Z} &= \eta Y - (\gamma+\alpha+\omega)Z+\alpha(N-1)[I],\\
\dot{\theta} &= -\frac{\beta[SI]}{Ng_x(\theta,t)}-(1-\theta)\left(\alpha\theta+\omega -\alpha(N-1)\frac{g(\theta,t)}{g_x(\theta,t)}\right),\\
\dot{\phi} &= 3\alpha - \left(\alpha + \omega+2\alpha(N-2)\frac{g_x(1,t)}{g_{xx}(1,t)}\right),\phi
\end{align}
where 
\begin{align}
    [SSI] &= [SS][SI]\frac{g_{xx}(\theta,t)}{N(g_x(\theta))^2}\left((1-\phi) + \phi g'(1)\frac{[SI]}{\theta g_x(\theta,t)Z}\right),\\
    [ESI] &= [SE][SI]\frac{g_{xx}(\theta,t)}{N(g_x(\theta,t))^2}\left((1-\phi) + \phi g'(1) N \frac{[EI]}{YZ}\right),\\
    [ISI] &= [SI]^2\frac{g_{xx}(\theta,t)}{N(g_x(\theta,t))^2}\left((1-\phi) + \phi g'(1) \frac{[II]}{\theta g_x(\theta,t)Y}\right).
\end{align}
%

\section{Additional Networks and Epidemiological\\ Parameters}\label{AppB}

In this appendix, we consider the prevalence-dependent intervention on two alternative heterogeneous, clustered networks with our COVID-19 parameters $R_0=2.4,\eta = 0.2$ and $\gamma = 0.1.$ We also consider the case where $R_0$ has increased to $5$ on the same networks, as well as the original unipartite projection contact network from the main text. The two networks considered are a Watts-Strogatz ``small world" network and a power law network with clustering. Both networks consist of $N=10,000$ nodes, as with the contact network in the main text. For the small world network, $\langle k \rangle \approx 30,\langle k^2-k\rangle \approx 900 ,\phi \approx 0.25$; for the the power law network with clustering, the relevant initial network parameters are $\langle k \rangle \approx 30,\langle k^2-k\rangle \approx 2000 ,\phi \approx 0.1.$ 

For each network, we report four sets of figures. First, we consider the relative change in final size (RCFS) and that average infections above threshold (AIAT) for the prevalence-dependent intervention with $q=0.01$ and $p=0.125,0.25,0.5.$ Each plot represents a combination of these two parameters, and $L_I$ and $L_R$ vary on the axes from $2$ to $180.$ Second, we fix $L_R=90$ and $L_I = 15, 30, 60$ and allow $p$ and $q$ to vary on the axes, with $p$ ranging from $0$ to $1$ and $q$ ranging from $0$ to $0.03$.

\begin{figure}[h!]
    \centering
    \begin{subfigure}{\textwidth}
        \includegraphics[width=\textwidth]{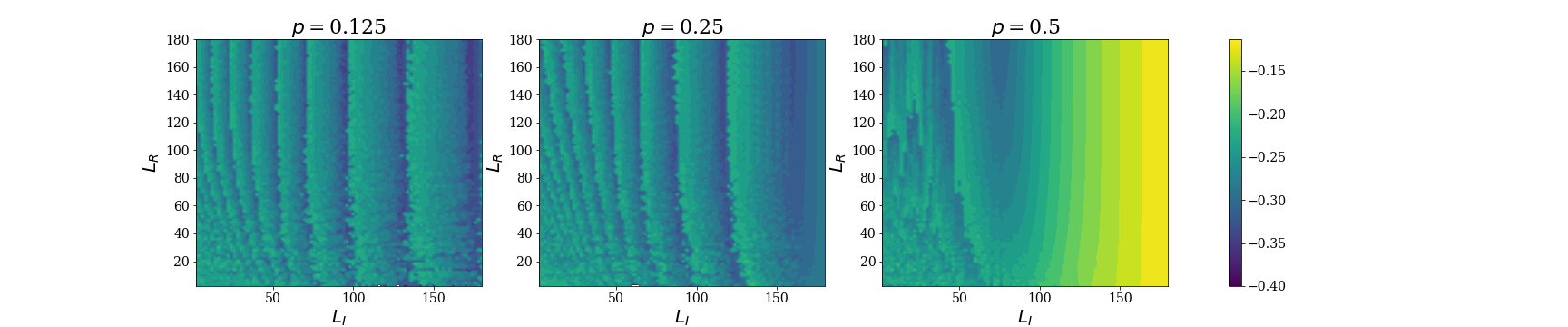}
        \caption{}
    \end{subfigure}
    \begin{subfigure}{\textwidth}
        \includegraphics[width=\textwidth]{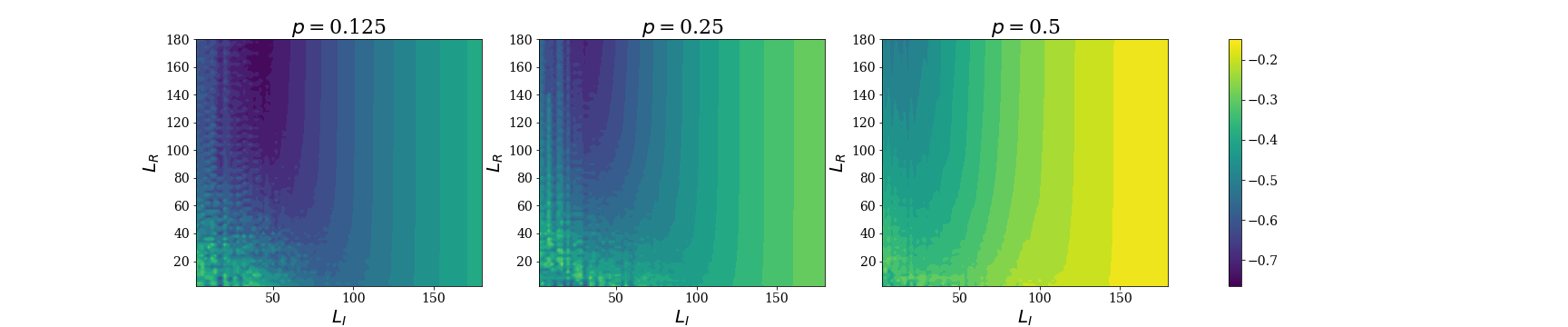}
        \caption{}
    \end{subfigure}
    \caption{Plots of the RCFS for the prevalence-dependent intervention on the (a) small world network and (b) power law network with clustering for $R_0 = 2.4,\eta=0.2,\gamma=0.1$ and a fixed $q=0.01$.}
    \label{fig:A1}
\end{figure}

\begin{figure}[h!]
    \centering
    \begin{subfigure}{\textwidth}
        \includegraphics[width=\textwidth]{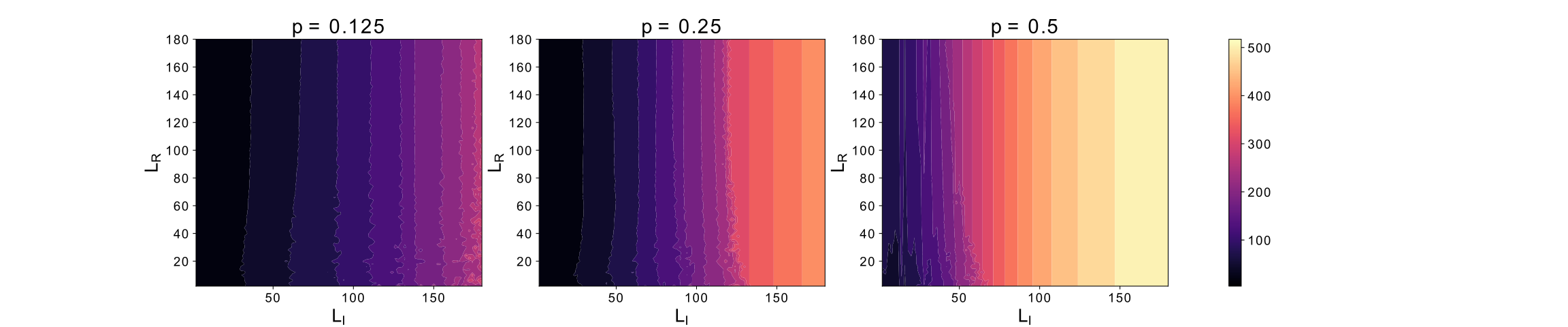}
        \caption{}
    \end{subfigure}
    \begin{subfigure}{\textwidth}
        \includegraphics[width=\textwidth]{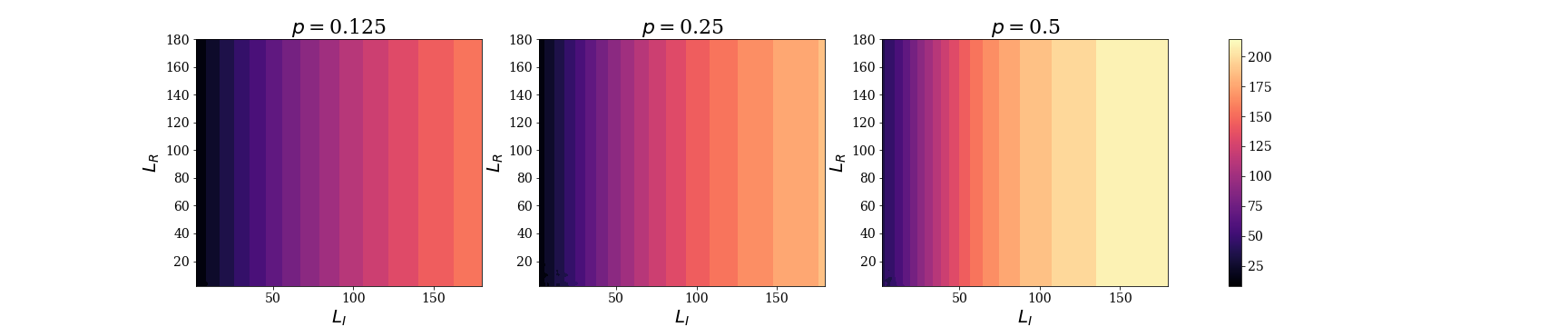}
        \caption{}
    \end{subfigure}
    \caption{Plots of the AIAT for the prevalence-dependent intervention on the (a) small world network and (b) power law network with clustering for $R_0 = 2.4,\eta=0.2,\gamma=0.1$ and a fixed $q=0.01$.}
    \label{fig:A2}
\end{figure}

\begin{figure}[h!]
    \centering
    \begin{subfigure}{\textwidth}
        \includegraphics[width=\textwidth]{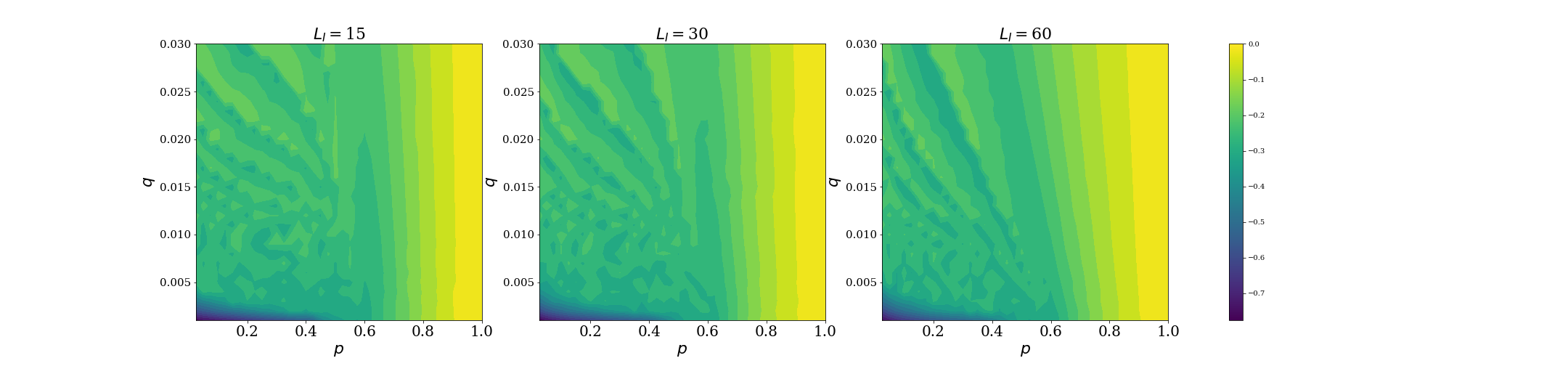}
        \caption{}
    \end{subfigure}
    \begin{subfigure}{\textwidth}
        \includegraphics[width=\textwidth]{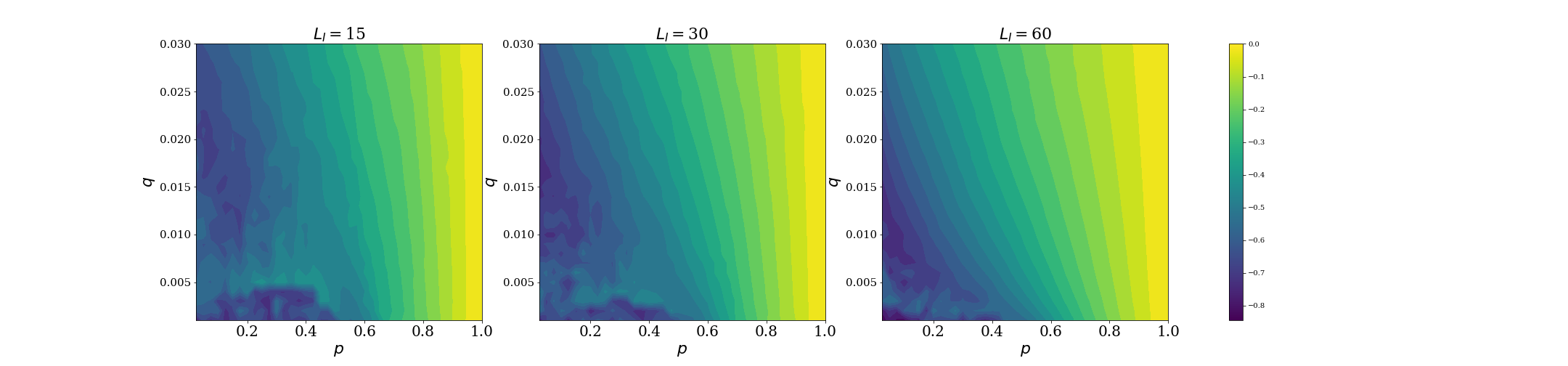}
        \caption{}
    \end{subfigure}
    \caption{Plots of the RCFS for the prevalence-dependent intervention on the (a) small world network and (b) power law network with clustering for $R_0 = 2.4,\eta=0.2,\gamma=0.1$ and a fixed $L_R=90$.}
    \label{fig:A3}
\end{figure}

\begin{figure}[h!]
    \centering
    \begin{subfigure}{\textwidth}
        \includegraphics[width=\textwidth]{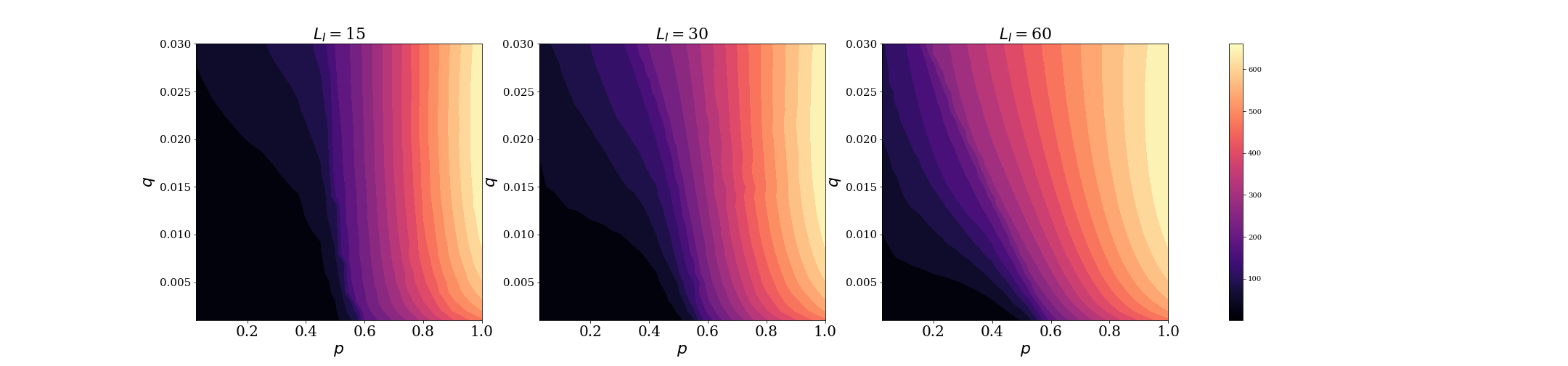}
        \caption{}
    \end{subfigure}
    \begin{subfigure}{\textwidth}
        \includegraphics[width=\textwidth]{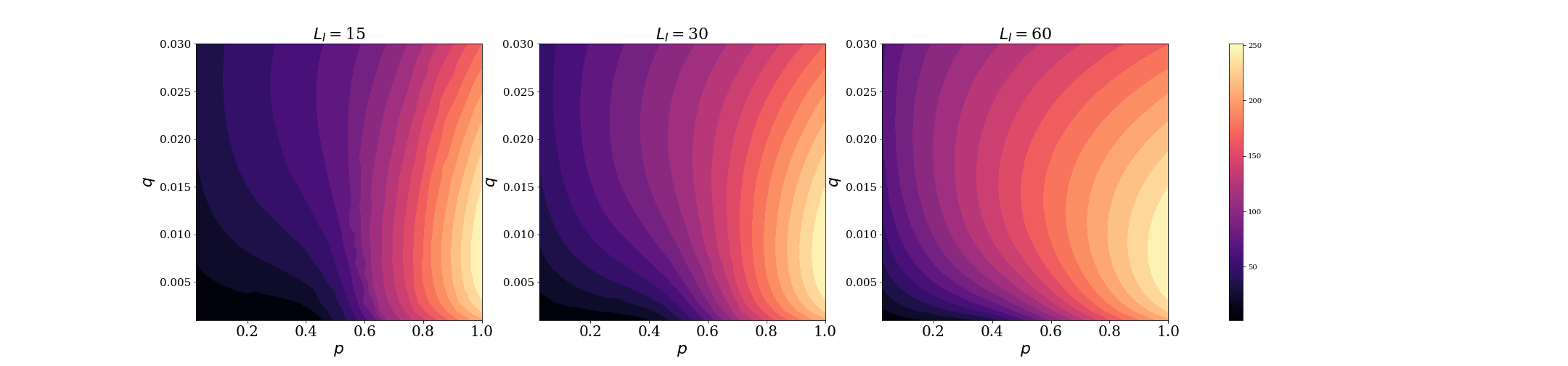}
        \caption{}
    \end{subfigure}
    \caption{Plots of the AIAT for the prevalence-dependent intervention on the (a) small world network and (b) power law network with clustering for $R_0 = 2.4,\eta=0.2,\gamma=0.1$ and a fixed $L_R = 90$.}
    \label{fig:A4}
\end{figure}

\begin{figure}
    \centering
    \begin{subfigure}{\textwidth}
    \centering
        \includegraphics[width=\textwidth]{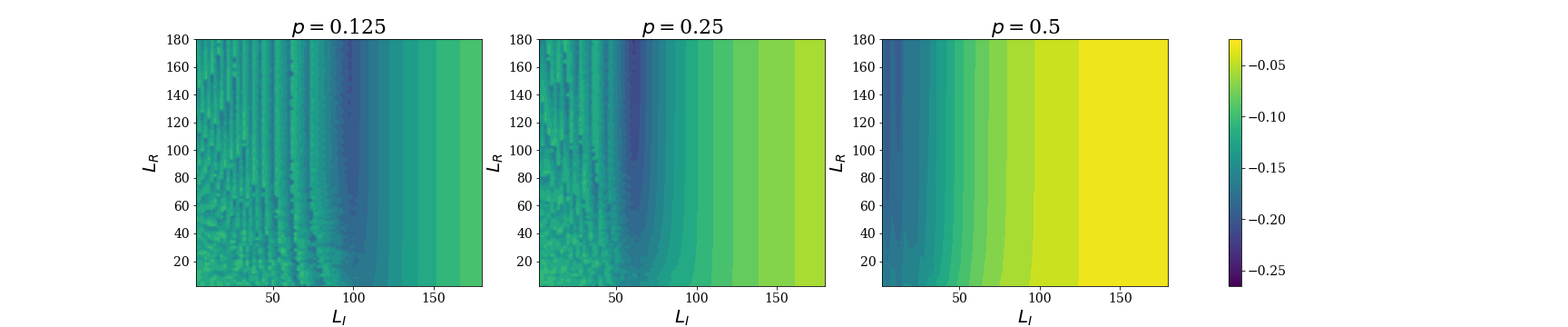}
        \caption{}
    \end{subfigure}
    \begin{subfigure}{\textwidth}
    \centering
        \includegraphics[width=\textwidth]{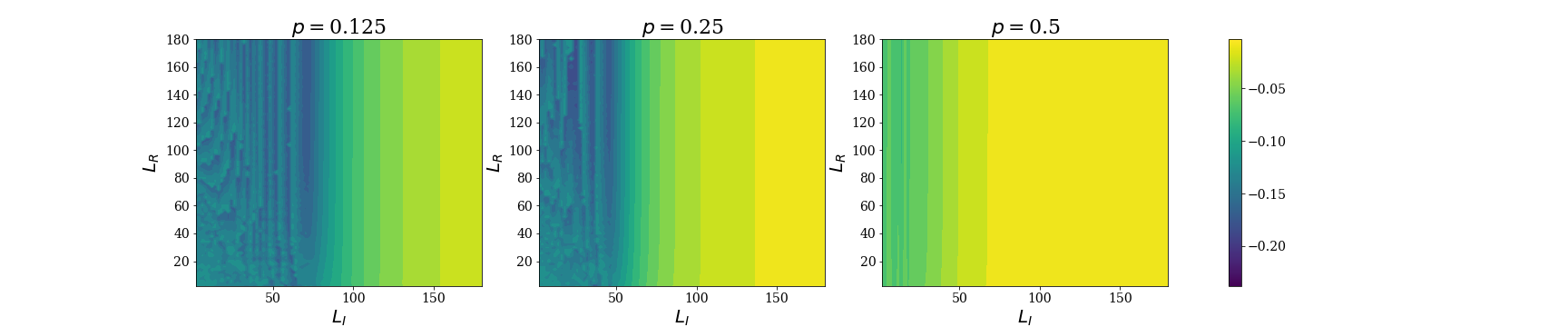}
        \caption{}
    \end{subfigure}
    \begin{subfigure}{\textwidth}
    \centering
        \includegraphics[width=\textwidth]{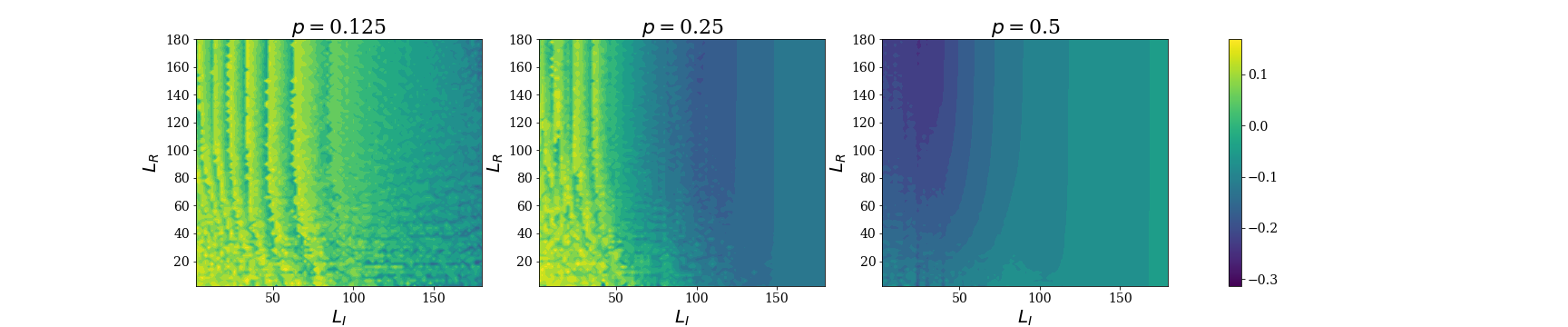}
        \caption{}
    \end{subfigure}
    \caption{Plots of the RCFS for the prevalence-dependent intervention on the (a) unipartite projection network, (b) small world network and (c) power law network with clustering for $R_0 = 5,\eta=0.2,\gamma=0.1$ and a fixed $q=0.01$.}
    \label{fig:A5}
\end{figure}

\begin{figure}
    \centering
    \begin{subfigure}{\textwidth}
    \centering
        \includegraphics[width=\textwidth]{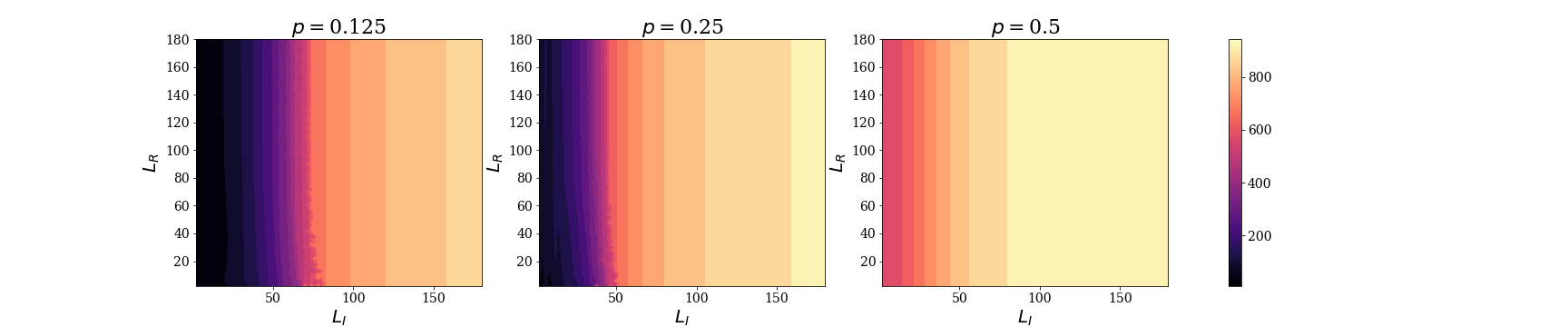}
        \caption{}
    \end{subfigure}
    \begin{subfigure}{\textwidth}
    \centering
        \includegraphics[width=\textwidth]{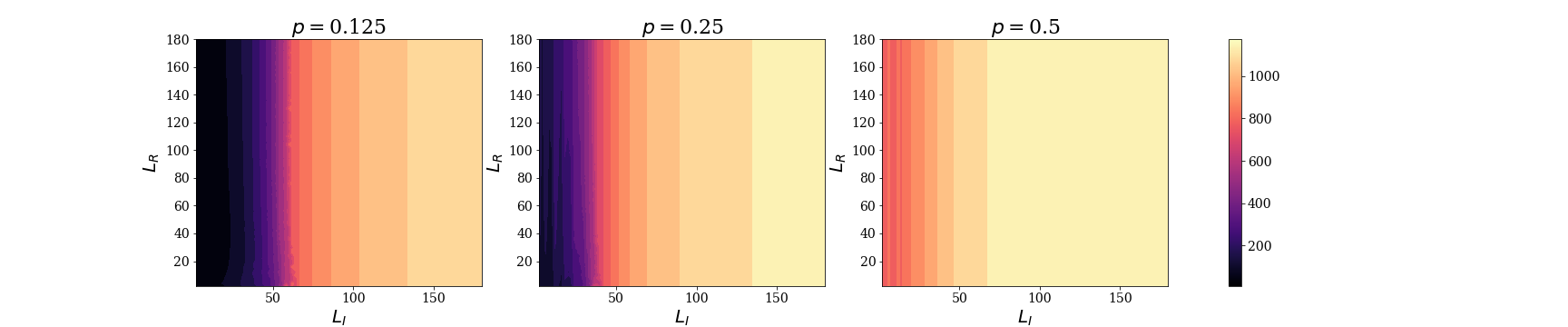}
        \caption{}
    \end{subfigure}
    \begin{subfigure}{\textwidth}
    \centering
        \includegraphics[width=\textwidth]{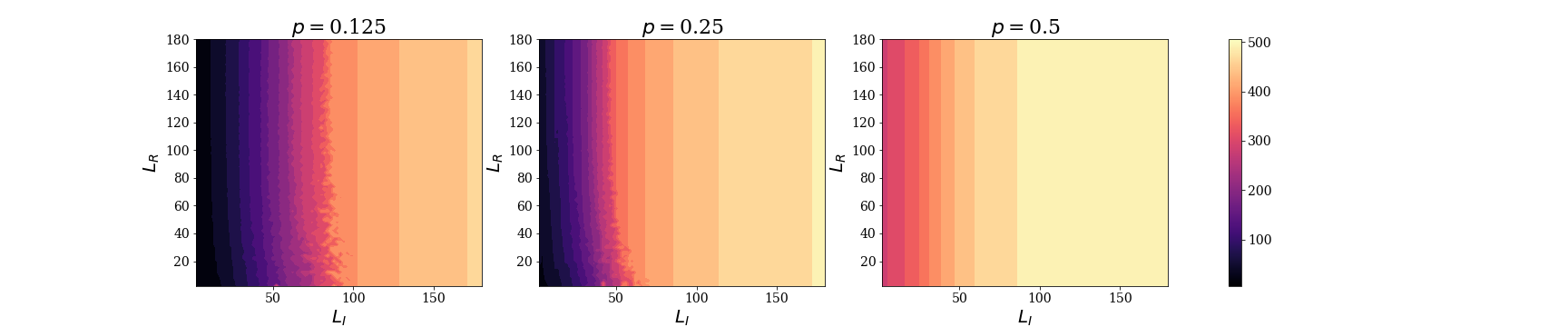}
        \caption{}
    \end{subfigure}
    \caption{Plots of the AIAT for the prevalence-dependent intervention on the (a) unipartite projection network, (b) small world network and (c) power law network with clustering for $R_0 = 5,\eta=0.2,\gamma=0.1$ and a fixed $q=0.01$.}
    \label{fig:A6}
\end{figure}

\begin{figure}
    \centering
    \begin{subfigure}{\textwidth}
    \centering
        \includegraphics[width=\textwidth]{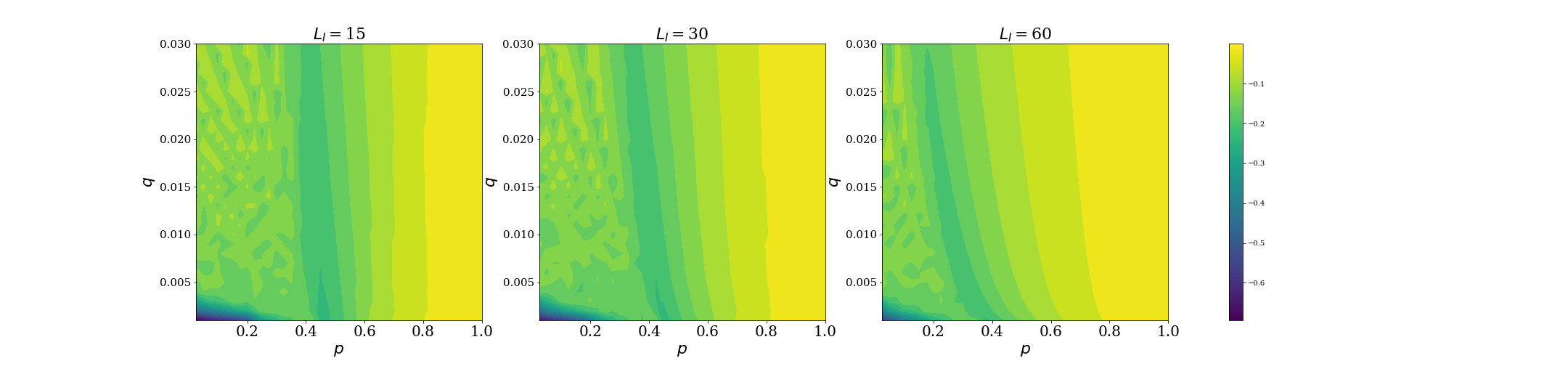}
        \caption{}
    \end{subfigure}
    \begin{subfigure}{\textwidth}
    \centering
        \includegraphics[width=\textwidth]{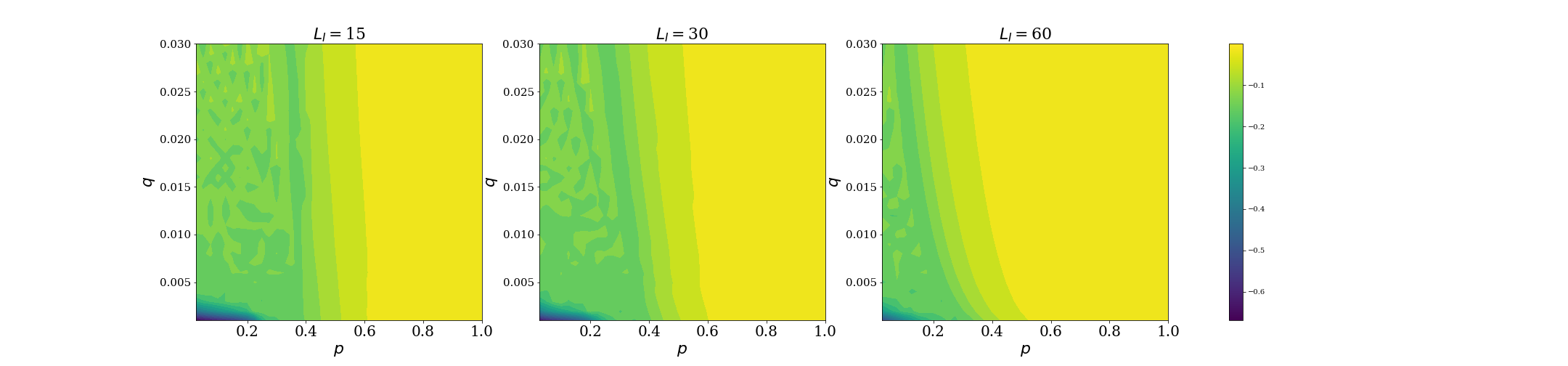}
        \caption{}
    \end{subfigure}
    \begin{subfigure}{\textwidth}
    \centering
        \includegraphics[width=\textwidth]{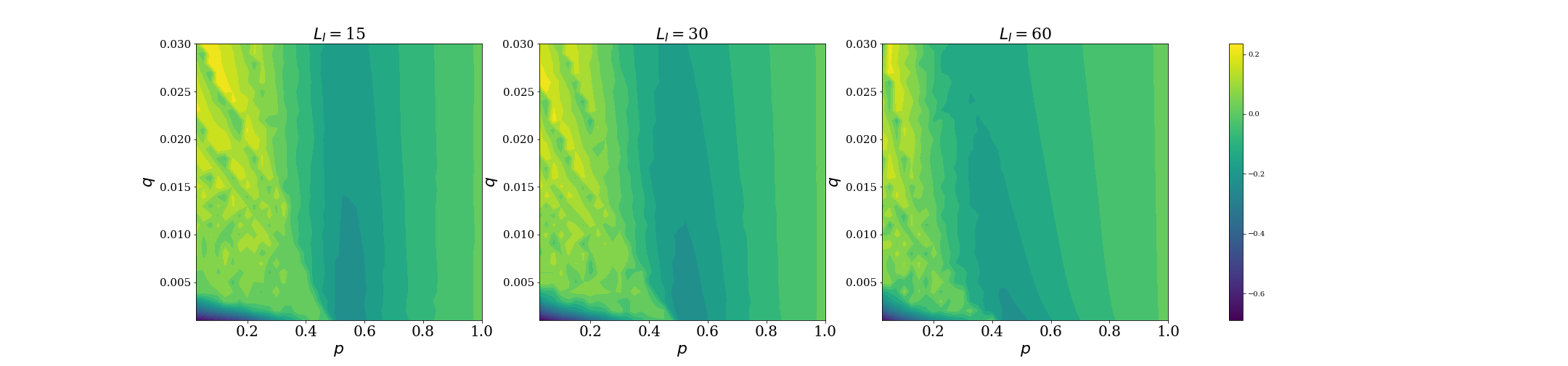}
        \caption{}
    \end{subfigure}
    \caption{Plots of the RCFS for the prevalence-dependent intervention on the (a) unipartite projection network, (b) small world network and (c) power law network with clustering for $R_0 = 5,\eta=0.2,\gamma=0.1$ and a fixed $L_R=90$.}
    \label{fig:A7}
\end{figure}

\begin{figure}
    \centering
    \begin{subfigure}{\textwidth}
    \centering
        \includegraphics[width=\textwidth]{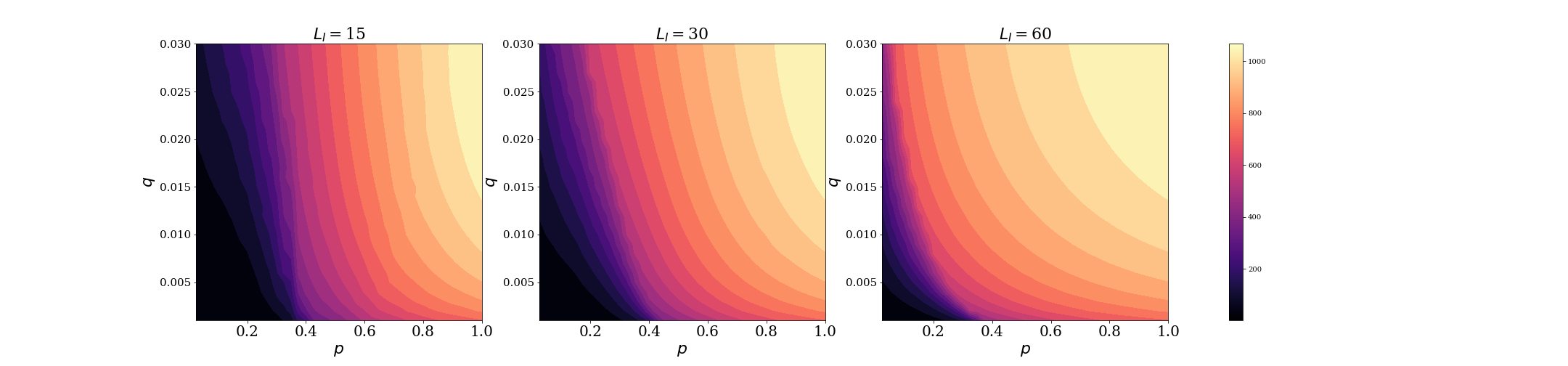}
        \caption{}
    \end{subfigure}
    \begin{subfigure}{\textwidth}
    \centering
        \includegraphics[width=\textwidth]{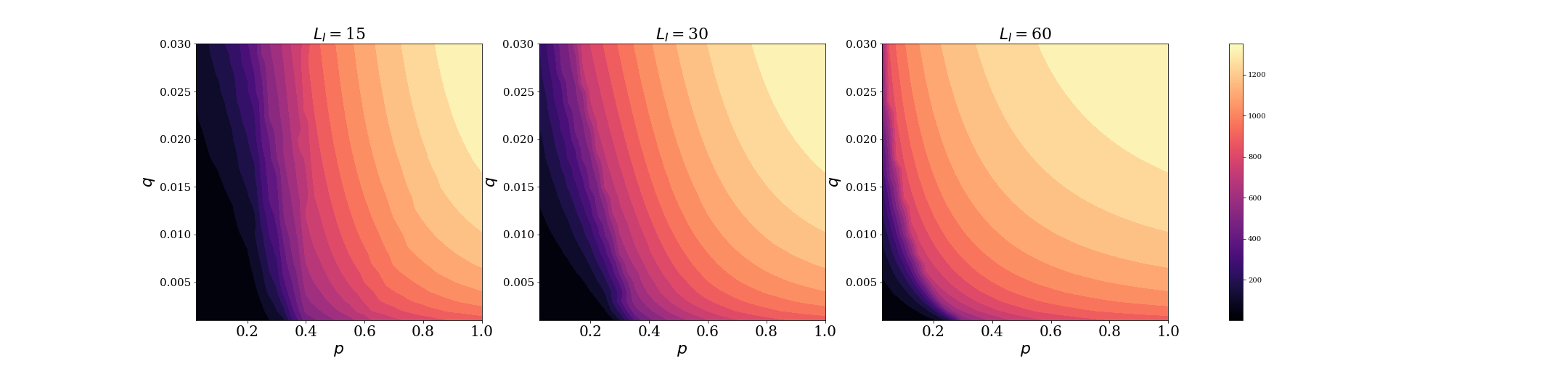}
        \caption{}
    \end{subfigure}
    \begin{subfigure}{\textwidth}
    \centering
        \includegraphics[width=\textwidth]{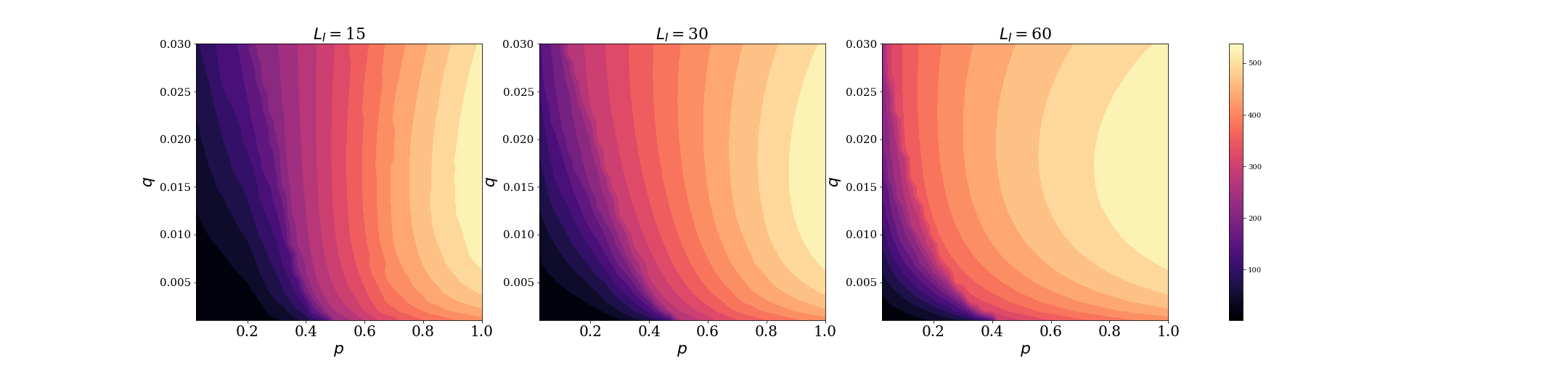}
        \caption{}
    \end{subfigure}
    \caption{Plots of the AIAT for the prevalence-dependent intervention on the (a) unipartite projection network, (b) small world network and (c) power law network with clustering for $R_0 = 5,\eta=0.2,\gamma=0.1$ and a fixed $L_R=90$.}
    \label{fig:A8}
\end{figure}

\end{appendices}

\end{document}